\documentclass[12pt]{amsart}

\input xy
\xyoption{all}

\usepackage{geometry}
\usepackage{amsmath}
\usepackage{amssymb}
\usepackage{amsthm}  
\usepackage{cite}
\newtheorem{same}{This should never appear}[section]
\newtheorem{defin}[same]{Definition}

\newtheorem{remark}[same]{Remark}
\newtheorem{theorem}[same]{Theorem}
\newtheorem{example}[same]{Example}
\newtheorem{lemma}[same]{Lemma}
\newtheorem{fact}[same]{Fact}

\newtheorem{cor}[same]{Corollary}
\newtheorem{prop}[same]{Proposition}

\newtheorem{conj}[same]{Conjecture}

\newbox\noforkbox \newdimen\forklinewidth
\setbox0\hbox{$\textstyle\smile$}
\setbox1\hbox to \wd0{\hfil\vrule width \forklinewidth depth-2pt
 height 10pt \hfil}
\setbox\noforkbox\hbox{\lower 2pt\box1\lower 2pt\box0\relax}
\def\unionstick{\mathop{\copy\noforkbox}\limits}

\def\nonfork_#1{\unionstick_{\textstyle #1}}

\setbox0\hbox{$\textstyle\smile$}
\setbox1\hbox to \wd0{\hfil{\sl /\/}\hfil}
\setbox2\hbox to \wd0{\hfil\vrule height 10pt depth -2pt width
              \forklinewidth\hfil}
\newbox\doesforkbox
\setbox\doesforkbox\hbox{\lower 2pt\box1 \lower 2pt\box2\lower2pt\box0\relax}
\def\nunionstick{\mathop{\copy\doesforkbox}\limits}

\def\fork_#1{\nunionstick_{\textstyle #1}}

\newcommand{\im}{\te{im }}

\newcommand{\sea}{\mathfrak{C}}

\newcommand{\prf}{\textbf{Proof: }}
\newcommand{\ba}{\bold{a}}

\newcommand{\emm}{\mathcal{M}}
\newcommand{\enn}{\mathcal{N}}
\newcommand{\eff}{\mathcal{F}}

\newcommand{\bN}{\bold{N}}

\newcommand{\dom}{\textrm{dom }}

\newcommand{\cf}{\text{cf }}
\newcommand{\rest}{\upharpoonright}

\newcommand{\seq}[1]{\langle #1 \rangle}
\newcommand{\te}[1]{\textrm{#1}}

\title{Tameness from Large Cardinal Axioms}
\author{Will Boney}
\email{wboney@cmu.edu}
\address{Department of Mathematical Sciences \\ Carnegie Mellon University \\ Pittsburgh, Pennsylvania, USA}

\date{May 8, 2014} 

\begin{document}

\maketitle

\begin{abstract}
We show that Shelah's Eventual Categoricity Conjecture for successors follows from the existence of class many strongly compact cardinals.  This is the first time the consistency of this conjecture has been proven.  We do so by showing that every AEC with $LS(K)$ below a strongly compact cardinal $\kappa$ is $< \kappa$-tame and applying the categoricity transfer of Grossberg and VanDieren \cite{tamenessthree}.  These techniques also apply to measurable and weakly compact cardinals and we prove similar tameness results under those hypotheses.  We isolate a dual property to tameness, called \emph{type shortness}, and show that it follows similarly from large cardinals.
\end{abstract}


\section{Introduction}

The study of Abstract Elementary Classes (AECs) began with Shelah's work in \cite{sh88} as a semantic generalization of the model theory of $L_{\lambda^+, \omega}(Q)$.  One of the main test questions for AECs is an attempt to prove an analogue of Morley's and Shelah's Categoricity Theorems \cite{morleycat} \cite{sh31} from first order logic to the AEC context.  This is typically referred to as Shelah's Categoricity Conjecture.  We state one of the more general versions, Shelah's Eventual Categoricity Conjecture that appears in the list of open problems from \cite{shelahfobook} as D.(3a).

\begin{conj}[Shelah]
For every $\lambda$, there is some $\mu_\lambda$ so that if $K$ is an AEC with $LS(K) = \lambda$ and is categorical in a cardinal greater than or equal to $\mu_\lambda$, then it is categorical in every cardinal greater than or equal to $\mu_\lambda$.
\end{conj}

Note that this is still open for countable fragments of $L_{\omega_1, \omega}$, where it is also conjectured that $\mu_{\aleph_0} = \beth_{\omega_1}$.  Shelah and others have made progress on this.  Some of this work also uses additional axioms of set theory, especially the work on frames contained in Shelah's recent book \cite{shelahaecbook} which uses many instances of the weak continuum hypothesis and non-saturation of the weak diamond ideals.  This work, using good $\lambda$-frames, is notable in that it focuses on transferring nice properties in small cardinals upwards with no global assumptions on the class.  Some work on the eventual categoricity problem has also been done from large cardinal axioms.  \cite{measure2} proves a downward categoricity transfer for inifinitary logics from the existence of a measureable cardinal and \cite{makkaishelah} proves Shelah's Categoricity Conjecture for $L_{\kappa, \omega}$ with $\kappa$ strongly compact.

Grossberg and VanDieren \cite{tamenessthree} have recently approached this problem from an exciting new approach.  They isolated a model theoretic property called tameness.  This is defined in Definition \ref{tamedef}, but briefly says that different types are different over small models.  Building on the results of \cite{sh394}, they showed that an upward version of Shelah's Categoricity Conjecture holds for tame AECs.

\begin{theorem}[\cite{tamenessthree}]\label{tamecatthm}
Suppose $K$ is an AEC with amalgamation, joint embedding, and no maximal models. If $K$ is $\chi$-tame and $\lambda^+$ categorical for some $\lambda \geq LS(K)^+ + \chi$, then $K$ is $\mu$ categorical for all $\mu \geq \lambda$.
\end{theorem}

An astute reader will notice the requirement that the categoricity cardinal be a successor cardinal that is not present in Morley's original theorem; this is a feature of most known cases of Shelah's Categoricity Conjecture and it is asked by Shelah in \cite{sh702}.6.14 if the successor requirement can be removed.  The results without the successor are Hytinnen and Kes\"{a}l\"{a}'s \cite{superstabfinitary} and \cite{cattransfinitary}, where they work with a strong assumptions of simplicity and finitarity.

This is the approach we use to prove the result stated in the introduction.  Briefly, the main theorem of this paper (Theorem \ref{strongcompactnesstheorem}) states

\begin{theorem}
If $K$ is an AEC with $LS(K) < \kappa$ and $\kappa$ is strongly compact, then $K$ is $\kappa$-tame.
\end{theorem}

This is combined with the result of Grossberg and VanDieren to give us the consistency of Shelah's Categoricity Conjecture.  Section \ref{conclusion} provides more details, including the derivation of amalgamation.  This result also improves \cite{makkaishelah} by showing that their categoricity result holds even for AECs that are not axiomatized by an infinitary theory.

The above result is part of a larger investigation of tame AECs that began with Grossberg and VanDieren's introduction of tameness in \cite{tamenessone}, which came from the latter's Ph.D. thesis.  In first-order model theory, types are trivially tame because different types necessarily contain different formulas with finitely many parameters, but, in an AEC, a type is not determined by formulas and instead have a semantic characterization (see Section \ref{prelimsection} below), so the question of tameness is not straightforward.  In the introduction of \cite{tamenessthree}, Grossberg and VanDieren list several previously studied nonelementary classes that turn out to be tame.  This list includes previous AECs for which a classification theory exists.  This lead them to the following conjecture about categoricity and tameness.

\begin{conj} \label{ramiconj}[Grossberg-VanDieren]
Suppose $K$ is an AEC.  If $K$ is categorical in some $\lambda \geq Hanf(LS(K))$ (or some other value depending only on $LS(K)$), then there exists $\chi < Hanf(LS(K))$ so that $K$ is $\chi$-tame.
\end{conj}

Our main theorem can be seen as proving a stronger version of this from the existence of a strongly compact cardinal instead of the categoricity assumption.

Some assumption (categoricity, large cardinals, etc.) is known to be necessary for any theorem that concludes tameness for many AECs.  This follows from the existence of nontame AECs.  Hart and Shelah \cite{hash323} implicitly provided the first example of a nontame AEC by constructing an infinitary sentence $\phi_k \in L_{\omega_1, \omega}$ for $k< \omega$ that is categorical at and before $\aleph_{k-2}$, but nowhere above $\aleph_{k-2}$.  Baldwin and Kolesnikov \cite{untame} later clarified this example by examining it specifically with tameness in mind and were able to show the exact failure of tameness.  Baldwin and Shelah \cite{nonlocality} created a more algebraic counterexample from short exact sequences of almost free, non-Whitehead groups of size $\kappa$.  Combining this with our result, this gives a new proof of the nonexistence of almost free, non-Whitehead groups above a strongly compact cardinal.  See Section \ref{furtherwork} for further discussion of these ideas.

In addition to tameness, we introduce and consider a dual locality property that we call \emph{type shortness}.  This property is defined explicitly in Section \ref{tametssection}, but briefly says that if two types of long, infinite sequences indexed by the same set differ, then there is a short subsequence where they already differ.  Comparing this with tameness, we are replacing the condition on the domain of the type with a condition on the index of the type realizers.  Type shortness and tameness together give a strong locality condition for wether a mapping can be extended to a $K$-embedding or an automorphism of the monster model.  We discuss this in depth in Section 3.  Additionally, the combination of these properties can be used to obtain a new notion of nonforking that can be seen as an AEC analogue of coheir.  This is done in \cite{shorttamedep}.

We now outline the paper.  Section 2 provides the AEC definitions and preliminaries that are necessary for this paper.  The only nonstandard item is that Galois types are allowed to be infinite in length (Definition \ref{galtype}) and Defintion \ref{essbelow} of ``essentially below,'' which captures exactly which AECs our results holds for.  Section 3 gives the various definitions of tameness and type shortness.  The main results of this paper are in Sections 4, 5, and 6.  Each of these sections assumes a different large cardinal axiom and uses a different technique to prove various levels of type shortness and tameness: Section 4 uses the ultrafilter definition of a strongly compact cardinal, Section 5 uses the elementary embedding definition of a measurable cardinal, and Section 6 uses the indescribability definition of a weakly compact cardinal.  Section 7 combines the results from this paper with the papers mentioned in the introduction.  This contains Theorem \ref{shelahcatconj}, the consistency of Shelah's Eventual Categoricity Conjecture for Successors.  Finally, Section 8 poses some new questions, especially in the area of the large cardinal strength of different universal tameness properties.

This paper was written while working on a Ph.D. under the direction of Rami Grossberg at Carnegie Mellon University and I would like to thank Professor Grossberg for his guidance and assistance in my research in general and in this work specifically.  I would also like to thank James Cummings, Ernest Schimmerling, and Spencer Unger for their discussions about the set theory and algebra involved in this paper; John Baldwin for reading an early version of this paper and pointing out an omission in Theorem \ref{bigresult}; the referee for many helpful comments; and my wife Emily Boney for her support.  After circulating a preprint, Shelah pointed us to his preprint \cite{sh932}, which contains similar ideas and proves some results independently.

\section{Preliminaries}

\label{prelimsection}

The definition for an Abstract Elementary Class was first given by Shelah in \cite{sh88}.  The definitions and concepts in the section are all part of the literature; in particular, see the books by Baldwin \cite{baldwinbook} and Shelah \cite{shelahaecbook}, the survey article by Grossberg \cite{grossberg2002}, or the forthcoming book by Grossberg \cite{ramibook} for general information

\begin{defin}
We say that $(K, \prec_K)$ is an Abstract Elementary Class iff
\begin{enumerate}

    \item There is some language $L = L(K)$ so that every element of $K$ is an $L$-structure; 

    \item $\prec_K$ is a partial order on $K$;

    \item for every $M, N \in K$, if $M \prec_K N$, then $M \subseteq N$;

    \item $(K, \prec_K)$ respects $L$ isomorphisms; that is, if $f: N \to N'$ is an $L$ isomorphism and $N \in K$, then $N' \in K$ and if we also have $M \in K$ with $M \prec_K N$, then $f(M) \in K$ and $f(M) \prec_K N'$;

    \item \emph{(Coherence)} if $M_0, M_1, M_2 \in K$ with $M_0 \prec_K M_2$, $M_1 \prec_K M_2$, and $M_0 \subseteq M_1$, then $M_0 \prec M_1$;

    \item \emph{(Tarski-Vaught axioms)} suppose $\seq{M_i \in K : i < \alpha}$ is a $\prec_K$-increasing continuous chain, then

        \begin{enumerate}

            \item $\cup_{i < \alpha} M_i \in K$ and, for all $i < \alpha$, we have $M_i \prec_K \cup_{i < \alpha} M_i$; and

            \item if there is some $N \in K$ so that, for all $i < \alpha$, we have $M_i \prec_K N$, then we also have $\cup_{i < \alpha} M_i \prec_K N$.; and

        \end{enumerate}

    \item \emph{(Lowenheim-Skolem number)} $LS(K)$ is the minimal cardinal $\lambda \geq |L(K)| + \aleph_0$ such that for any $M \in K$ and $A \subset |M|$, there is some $N \prec_K M$ such that $A \subset |N|$ and $\|N\| \leq |A| + \lambda$.

\end{enumerate}

\end{defin}

\begin{remark}
 As is typical, we drop the subscript on $\prec_K$ when it is clear from context and we abuse notation by calling $K$ an AEC when we mean $(K, \prec_K)$ is an AEC.  Also, we follow the convention of Shelah that, for $M \in K$, we differentiate between the model $M$, its universe $|M|$, and the cardinality of its universe $\|M\|$.  In this paper, $K$ is always an AEC and that has no models of size smaller than the Lowenheim-Skolem number.
\end{remark}

The class of AECs is strong enough to encompass commonly studied logical frameworks, including classes  of models of theories in first-order logc, infinitary logics, logics with added quantifiers, and more.

After seeing a preliminary version of this work, Jose Iovino pointed us to the work on Metric Abstract Elementary Classes (MAECs) by Hirvonen and Hyttinen \cite{catmaecs} and others.  This is a more general framework that extends AECs as continuous first-order logic extends first-order logic and is more suited for dealing with analytic concepts like being a complete metric space.  Although there is a slightly different notion of ultraproducts, the theorems of this paper still hold in that context.

We will briefly summarize some of the basic notations, definitions, and results for AECs; as above, see \cite{ramibook} for a more detailed description and development.

\begin{defin}
\begin{enumerate}

	\item A \emph{$K$ embedding} from $M$ to $N$ is an injective map $f: M \to N$ so $f(M) \prec_K N$.

	\item \begin{eqnarray*}
	K_\lambda &=& \{ M \in K : \|M\| = \lambda \}\\
	K_{\leq \lambda} &=& \{ M \in K : \|M\| \leq \lambda \}
	\end{eqnarray*}

	\item $K$ has the \emph{$\lambda$-amalgamation property} ($\lambda$-AP) iff for any $M \prec N_0, N_1 \in K_\lambda$, there is some $N^* \in K$ and $f_i: M \to N_i$ so that
\[
 \xymatrix{\ar @{} [dr] N_1  \ar[r]^{f_1}  & N^*\\
M \ar[u] \ar[r] & N_0 \ar[u]_{f_2}
 }
\]
commutes.
	
	\item $K$ has the \emph{$\lambda$-joint embedding property} ($\lambda$-JP) iff for any $M_0, M_1 \in K_\lambda$, there is some $M^* \in K$ and $f_i: M_i \to M^*$.  This is also called the joint embedding property.
	
	\item $K$ has \emph{no maximal models} iff for every $M \in K$, there is some $N \in K$ so $M \precneqq N$.

	\item Fix $\lambda$.  Set 
	\begin{eqnarray*}
	Hanf(\lambda) &=& \{ \mu^+ : \te{ there is an AEC with $LS(K)$ that has a model of size $\mu$ }\\
	& & \te{but not arbitrarily large models }\}
	\end{eqnarray*}
	We sometimes write $Hanf(K)$ for $Hanf(LS(K))$.

\end{enumerate}
\end{defin}

In AECs, types as sets of formulas do not behave as nicely as they do in first-order model theory; any of the examples of non-tameness is an example of this and it is made explicit in \cite{untame}.  However, Shelah isolated a semantic notion of type in \cite{sh300} that Grossberg named \emph{Galois type} in \cite{grossberg2002} this can replace the first-order notion.

We differ from the standard treatment of types in that we allow the length of our types to be possibly infinite.  This is useful because it is sometimes natural, as in \cite{makkaishelah},\cite{shelahaecbook}.V, or \cite{longtypes}, and because it allows us to show the full power of our results in the following sections.  For early results in the stability theory of $\alpha$-types, see \cite{tamenessone}.

\begin{defin} \label{galtype} Let $K$ be an AEC, $\lambda \geq LS(K)$, and $(I, <_I)$ an ordered set.
 \begin{enumerate}
  
    \item Set $K^{3, I}_\lambda = \{ (\seq{a_i : i \in I}, M, N) : M \in K_\lambda, M \prec N \in K_{\lambda + |I|}, \te{ and } \{a_i : i \in I \} \subset |N| \}$.  The elements of this set are called \emph{pretypes}.

    \item Given two pretypes, $(\seq{a_i : i \in I}, M, N)$ and $(\seq{b_i : i \in I}, M', N') \in K_\lambda^{3, I}$, we say that $(\seq{a_i : i \in I}, M, N) \sim_{AT} (\seq{b_i : i \in I}, M', N')$ iff $M = M'$ and there is $N^* \in K$ and $f: N \to N^*$ and $g: N' \to N^*$ so that $f(a_i) = g(b_i)$ for all $i \in I$ and the following diagram commutes:
\[
 \xymatrix{\ar @{} [dr] N'  \ar[r]^{g}  & N^*\\
M \ar[u] \ar[r] & N \ar[u]_{f}
 }
\]

    \item Let $\sim$ be the transitive closure of $\sim_{AT}$.
    
    \item For $M \in K$, set $tp(\seq{a_i : i \in I}/M, N) = [(\seq{a_i : i \in I}, M, N)]_\sim$ and $S^I(M) = \{ tp(\seq{a_i : i \in I}/M, N) : (\seq{a_i : i \in I}/M, N) \in K^{3, I}_{\| M \|} \}$.
    
    \item Suppose $p = tp(\seq{a_i : i \in I}/M, N)  \in S^I(M)$ and $M' \prec M$ and $I' \subset I$.  Then, 
    \begin{itemize}
    
    	\item $p \rest M' \in S^I(M')$ is $tp(\seq{a_i : i \in I}/M', N')$ for some (any) $N' \in K_{\|M'\|+ |I|}$ with $M' \prec N' \prec N$ and $\seq{a_i : i \in I} \subset |N'|$;
    	
    	\item and $p^{I'} \in S^{I'}(M)$ is $ tp(\seq{a_i : i \in I'}/M, N)$.
 
	\end{itemize} 
    
\end{enumerate}

\end{defin}

\begin{remark}
\begin{enumerate}

	\item If $K$ has the $\lambda + |I|$-amalgamation property, then $\sim_{AT}$ is a transitive relation and, thus, an equivalence relation on $K^{3, I}_\lambda$; note that `$AT$' stands for ``atomic.''
	
	\item Some authors place a $g$ denoting `Galois' in front of the above notions to differentiate them from the first-order versions (ie, $gtp(a/M, N)$ and $gS(M)$); however, since we almost exclusive use  Galois types and only reference syntactic types in this section, we omit this.

\end{enumerate}
\end{remark}

In the presence of sufficiently strong AP and JP and if $K$ has no maximal models, there exists a monster model, which greatly simplifies the notion of type.

\begin{defin}[Half definition/half remark]
Suppose $K$ is an AEC with amalgamation, joint embedding, and no maximal models.  There is a monster model $\sea$; that is, a model from $K$ of large size of high cofinality that is universal and model homogeneous for all models that we will consider; that is all $N \in K$ can be embedded into $\sea$ and if $M \prec \sea$ and $M \prec N$, then there is some $f: N \underset{M}{\to} \sea$.
\end{defin}

See \cite{ramibook}.4.4 for a more detailed discussion of monster models.

\begin{defin}
If $K$ has a monster model $\sea$, then $gtp(a/M)$ is the orbit of $a$ under the action of automorphism of $\sea$ fixing $M$.  That is, $a$ and $b$ realize the same type over $M$ (equivalently, $gtp(a/M) = gtp(b/M)$) iff there is some $f \in Aut_M \sea$ so $f(a) = b$.
\end{defin}

These two notions of type are equivalent.  Note that this definition explains the name, as the orbits of automorphisms fixing smaller structures recalls certain aspects of Galois Theory.

Finally, we state Shelah's Presentation Theorem from \cite{sh88} that characterizes AECs as pseduoelementary classes and extends Chang's Presentation Theorem from \cite{chang68}.  This will be important for technical results later.

\begin{defin}
Let $T_1$ be a first-order theory, $\Gamma$ a set of finitary, syntactic $T_1$-types, and $L \subset L(T_1)$ a language.  The pseudoelementary class $PC(T_1, \Gamma, L) = \{ M \rest L : M \vDash T_1 \te{ and omits each } p \in \Gamma \}$ for a theory $T_1$, a set of $L(T_1)$ types $\Gamma$, and $L \subseteq L(T_1)$.  To say that $K$ is a $PC_{\lambda, \kappa}$ class means that $K = PC(T_1, \Gamma, L)$ for $|T_1| \leq \lambda$ and $|\Gamma| \leq \kappa$.
\end{defin}

\begin{theorem}[Shelah's Presentation Theorem]	Suppose $K$ is an AEC with $LS(K) = \kappa$.  There is some $L_1 \supseteq L(K)$ of size $\kappa$, a first-order theory $T_1$ in $L_1$ of size $\kappa$, and a set of $L_1$-types $\Gamma$ over the empty set (so $|\Gamma| \leq 2^\kappa$) so that $K = PC(T_1, \Gamma, L(K))$ and for any $M_1 \models T_1$ and $N_1 \subseteq M_1$, if $M_1$ omits $\Gamma$, then $N_1 \rest L(K) \prec_K M_1 \rest L(K)$.  Moreover, every $M \in K$ has an expansion to an $L(K)$ structure $M_1 \in EC(T_1, \Gamma)$ so that, for all $N$ that is an $L(K)$ structure,
	$$N \prec_K M \iff \te{ there is some $N_1 \subseteq M_1$ so $N = N_1 \rest L(K)$ }$$
\end{theorem}

We end the model theoretic preliminaries with a definition that will allow us to easily state which AECs are conclusions are valid for:

\begin{defin} \label{essbelow}
For a cardinal $\kappa$, we say that an AEC is \emph{essentially below $\kappa$} iff a) $LS(K) < \kappa$ or b) $K = (\te{Mod T}, \prec_\eff)$ for $T$ a $L_{\kappa, \omega}$ theory.
\end{defin}

We use heavily the standard ultraproduct construction.  Recall that, if $U$ is an ultrafilter on $I$, then
$$\Pi M_i / U = \{ [f]_U : f \in \Pi_{i \in I} M_i \}$$
where $[\cdot]_U$ denotes the equivalence class of a function under equality on a $U$-large set.

Finally, we recall two set-theoretic definitions used in applications of large cardinals that will be needed later.  For more detail, consult \cite{kanamori}.1.5.

\begin{defin}
\begin{itemize}
	\item If $j: V \to M$ is elementary with $M \neq V$ transitive, then $\text{crit } j = \min \{ \alpha \in ON : \alpha \neq j(\alpha)\}$; in fact, this is well-defined.
	
	\item An ultrafilter $U$ on $I$ is $\kappa$-complete iff, for any $\alpha < \kappa$ and $X_\beta \in U$ for  $\beta < \alpha$, we have $\cap_{\beta < \alpha} X_\beta \in U$.	
	
	\item An ultrafilter $U$ on $P_\kappa I$ (or a subset of it) is called \emph{fine} iff, for all $i \in I$, the set $[i] := \{ I_0 \in P_\kappa I : i \in I_0 \} \in U$.
	
	\item Given a well-founded set $(X, E)$, the \emph{Mostowski collapse} is unique function $\pi_X$ with domain $X$ so that $(\pi''X, \in)$ is transitive and, for all $y \in X$, $\pi(y) = \{ \pi(x) : x E y\}$.

\end{itemize}
\end{defin}

\section{Tameness and type shortness} \label{tametssection}

Tameness is a property first isolated by Rami Grossberg and Monica VanDieren in their papers \cite{tamenessone}, \cite{tamenesstwo}, and \cite{tamenessthree}; \cite{tamenesstwo} came from VanDieren's thesis \cite{monicathesis}.  The property is similar to one used by Shelah in \cite{sh394}, where he derived this property for types with \emph{saturated} domains from categoricity in a successor cardinal above the second Hanf number, $\beth_{(2^{\beth_{(2^{LS(K)})^+}})^+}$; this property is now called \emph{weak tameness} (see \cite{nonlocality}).  In their papers, Grossberg and VanDieren defined only $\chi$-tameness; the two cardinal parameterization of it appeared later in \cite{baldwinbook}.

We begin with a minor notational definition and then define several levels of tameness:

\begin{defin} \label{tamedef}
Suppose $K$ is an AEC with $LS(K) < \kappa \leq \lambda$ and $I$ is a linear order.
\begin{enumerate}

	\item For any $M \in K_{\geq \kappa}$, we write
	$$P^*_\kappa M = \{ N \prec M : \|N\| < \kappa\}$$

	\item $K$ is \emph{$(< \kappa, \lambda)$-tame for $I$-length types} iff for any $M \in K_\lambda$ and $p \neq q \in S^I(M)$, there is some $N \in P^*_{\kappa} M$ and $p \rest N \neq q \rest N$.
	
	\item $K$ is \emph{$< \kappa$-tame for $I$-length types} iff $K$ is $(< \kappa, \mu)$ tame for $I$-length types for all $\mu \geq \kappa$.

	\item $K$ is \emph{fully $<\kappa$-tame} iff $K$ is $<\kappa$-tame for $I$-length types for all $I$.

	\item Writing ``$\kappa$'' for ``$<\kappa$'' means ``$<\kappa^+$.''

\end{enumerate}	
\end{defin}

If we omit the $I$, we mean $I = 1$.  $P^*_\kappa M$ is reminiscent of the set theoretic notation $P_\kappa A = \{ X \subset A : |X| < \kappa \}$.

Note that we gave the above definitions as different types are different over a small model; this is clearly equivalent to saying that any two types which are the same over all small models are the same.  For instance,

\begin{fact}
$K$ is $< \kappa$-tame iff for any $M \in K_{\geq \kappa}$ and $p, q \in S(M)$, if $p \rest N = q \rest N$ for all $N \in P_\kappa^* M$, then $p = q$.
\end{fact}

Recall that, by definition, if the restrictions of two types to a smaller model are different, then the original types are different.  Tameness is a way of saying that the converse holds as well.
Obviously, if $K = (\te{Mod } T, \prec_{L_{\omega \omega}})$, then $LS(K) = |L(T)| + \aleph_0$ and $K$ is $LS(K)$ tame; in fact, given $p \neq q \in S(M)$, there is a finite tuple that witnesses their difference.

The power of tameness is shown through the categoricity transfer theorem of Grossberg and VanDieren from the introduction, Theorem \ref{tamecatthm}.

In this paper, we introduce a dual notion to tameness: \emph{type shortness}.  If we think of tameness as a locality property for the domains of types, then type shortness is a locality property for the \emph{length} of types.  Below we make this precise:

\begin{defin} Suppose $K$ is an AEC with $LS(K) \leq \mu$ and $\kappa < \lambda$.
\begin{enumerate}

	\item $K$ is \emph{$(< \kappa, \lambda)$-type short over $\mu$-sized models} iff for any $M \in K_\mu$ and $p \neq q \in S^\lambda(M)$, there is some $I' \subset I$ of size $< \kappa$ so that $p^{I'} \neq q^{I'}$.

	\item $K$ is \emph{$< \kappa$-type short over $\mu$-sized models } iff $K$ is $(< \kappa, \lambda)$-type short over $\mu$-sized models for all $\lambda \geq \kappa$.
	
	\item $K$ is \emph{fully $< \kappa$-type short} iff $K$ is $< \kappa$-type short over $\mu$-sized models for all $\mu$.
	
	\item Writing ``$\kappa$'' for ``$<\kappa$'' means ``$<\kappa^+$.''

\end{enumerate}	
\end{defin}

The reason for isolating type shortness is a bit more artificial than tameness: at the advice of Grossberg, we attempted to investigate an independence relation on tame classes following \cite{makkaishelah}.  In the course of doing so, this notion came to light.  Then, in revisiting the constructions in this paper, it was clear that they would provide large amounts of type shortness as well as tameness.  The results on the independence are described in Boney and Grossberg \cite{shorttamedep}.

The connection between tameness and type shortness is more than just a vague statement of duality.  Given varying strengths of one, we are able to get the other, as outlined in the two following theorems:

\begin{theorem}
Suppose $K$ is an AEC with amalgamation, joint embedding, and no maximal models.  If $K$ is categorical in $\mu$ and $(< \kappa, \mu)$-tame for $\lambda$-length types, then $K$ is $(< \kappa, \mu)$-type short for types of models over $\lambda$-sized domains.
\end{theorem}

{\bf Proof:} Let $M, M' \in K_\mu$ and $N \in K_\lambda$ so that $tp(M/N) \neq tp(M'/N)$.  By $\mu$ categoricity, there is some $f: M \cong M'$; WLOG $f \in Aut \sea$.\\
{\bf Claim:} $tp(f(N)/M') \neq tp(N/M')$\\
If not, then there is some $h \in Aut_{M'} \sea$ so $h \circ f(N) = N$.  Then $h \circ f \in Aut_N \sea$ and $h\circ f(M) = h(M') = M'$, which means $tp(M/N) = tp(M'/N)$, a contradiction. \hfill $\dag_{\te{Claim}}$\\
Now, by tameness, there is some $M^- \in P_\kappa^* M'$ so $tp(f(N)/M^-) \neq tp(N/M^-)$.  Then, by the same argument as in the claim, we get that $tp(f^{-1}(M^-)/N) \neq tp(M^-/N)$, which is what we want because $f^{-1}(M^-) \in P_\kappa^* M$. \hfill \dag

\begin{theorem}
Suppose $K$ is an AEC with amalgamation, joint embedding, and no maximal models.  If $K$ is $(< \kappa, \mu)$-type short over the empty set, then it is $(< \kappa, \mu)$-tame for $\leq \mu$-length types.
\end{theorem}

{\bf Proof:} Suppose $tp(a/M) \neq tp(b/M)$ for $M \in K_\mu$ and $\ell(a) = \ell(b) \leq \mu$.  Then we have $tp(aM/\emptyset) \neq tp(bM/\emptyset)$.  By our type shortness, there is some $a' \subset a$, $b' \subset b$, and $X_0 \subset M$ all of cardinality $< \kappa$ so that $tp(a' X_0/\emptyset) \neq tp(b' X_0/\emptyset)$.  Then find $M_0 \prec M$ of size $< \kappa$ that contains $X$.  Then
\begin{eqnarray*}
tp(a' X_0/\emptyset) &\neq& tp(b' X_0/\emptyset) \\
tp(a M_0/\emptyset) &\neq& tp(b M_0/\emptyset) \\
tp(a/M_0) &\neq& tp(b/M_0)
\end{eqnarray*}
as desired. \hfill \dag\\

The hypothesis that an AEC is both $<\kappa$ tame and $<\kappa$ type short, as in the conclusion of Theorem \ref{strongcompactnesstheorem}, gives a locality condition for testing wether a function $f$ that fixes a model $M$ is or can be extended to a $K$ embedding.  If there is no $K$ embedding that fixes $M$ and sends $\dom f$ to $\im f$, then $tp(\dom f/M) \neq tp( \im f/M)$.  Then, by tameness and type shortness, there is some $M_0 \in P_\kappa^* M$ and $X_0 \in P_\kappa \dom f$ so we have $tp(X_0/M_0) \neq tp(f(X_0)/M_0)$.
Thus, in a $<\kappa$ tame and type short AEC, $f$ can be extended to a $K$ embedding iff every subset of $f$ of size $<\kappa$ can be extended to a $K$ embedding.  We explore these AECs more in \cite{shorttamedep}.

There are other properties of AECs that assert different locality properties of types.  \cite{nonlocality} contains some of these.  The arguments in the following sections are also useful in deriving those properties.

\section{Strongly Compact}
\label{stronglycompactsection}

We begin with a study of AECs under the assumptions that there is a strongly compact cardinal $\kappa$ and a given AEC is essentially below $\kappa$ (see Definition \ref{essbelow}), but has a model above $\kappa$.  Since $\kappa$ is strongly inaccessible, this is equivalent to the AEC having a model above its Hanf number.

\begin{defin}[\cite{jech}.20] \label{strongcompdef}
An uncountable cardinal $\kappa$ is \emph{strongly compact} iff every $\kappa$-complete filter can be extended to a $\kappa$-complete ultrafilter.\\
	Equivalently, $L_{\kappa, \omega}$ and $L_{\kappa, \kappa}$ satisfy the compactness theorem.\\
	Equivalently, for every $\lambda \geq \kappa$, there is some elementary (in the first-order sense) embedding $j: V \to \emm$ with critical point $\kappa$ such that $j(\kappa) > \lambda$ and there is some $Y \in \emm$ of size $\lambda$ such that $j''\lambda \subset Y$.\\
	Equivalently, for every $\lambda \geq \kappa$, there is a fine, $\kappa$ complete ultrafilter $U$ on $P_\kappa \lambda$; that is, a $\kappa$ complete ultrafilter so, for every $\alpha < \kappa$, we have  $[\alpha] = \{ X \in P_\kappa \lambda : \alpha \in X \} \in U$.
\end{defin}

In this section, we prefer to use the latter ultrafilter formulation because it is more model-theoretic in nature.  In the next section, on measurable cardinals, we discuss the elementary embedding formulation of a large cardinal that is preferred by set theorists.

The most basic and fundamental model-theoretic fact about ultraproducts is \L o\'{s}' Theorem, which tells us that $Mod \te{ } T$ is closed under ultraproducts.  We wish to prove a version of \L o\'{s}' Theorem for AECs.  This will generalize the version for $L_{\kappa, \omega}$ when $\kappa$ is measureable.

\begin{theorem}[\L o\'{s}' Theorem for $L_{\kappa, \omega}$]
Let $U$ be a $\kappa$-complete ultrafilter over $I$, $L$ be a language, and $\seq{M_i : i \in I}$ be $L$ structures.  Then, for any $[f_1]_U, \dots [f_n]_U \in \Pi M_i / U$ and $\phi(x_1, \dots, x_n) \in L_{\kappa, \omega}$, we have
$$\Pi M_i / U \models \phi([f_1]_U, \dots, [f_n]_U) \te{ iff } \{ i \in I : M_i \models \phi(f_1(i), \dots, f_n(i)) \} \in U$$
\end{theorem}

This is proved similarly to the first-order version.  Our version for AECs is necessarily more complex since we do not have any syntax.  Thus, the characterization must be done semantically.  However, the following theorem aims to obtain the same results as the first-order version.  Of particular interest are parts (5) and (6): (5) says that if $M = \cup_{i < \kappa} M_i$ for $\seq{M_i : i < \kappa}$ increasing, then we can cannonically embed $M$ into $\Pi M_i / U$ and (6) says the same thing for $\seq{M_i : i < \kappa}$ a directed set.

\begin{theorem}[\L o\'{s}' Theorem for AECs] \label{losaec}
Suppose $K$ is an AEC essentially below $\kappa$ and $U$ is a $\kappa$-complete ultrafilter on $I$.  Then $K$ and the class of $K$-embeddings is closed under $\kappa$-complete ultrapowers and the ultrapower embedding.  In particular,
\begin{enumerate}

	\item if $\seq{M_i \in K : i \in I}$, then $\Pi M_i / U \in K$;
	
	\item if $\seq{M_i \in K : i \in I}, \seq{N_i \in K : i \in I}$ and, for every $i \in I$, $M_i \prec_K N_i$, then $\Pi M_i / U \prec_K \Pi N_i / U$
	
	\item if $\seq{M_i \in K : i \in }, \seq{N_i \in K : i \in I}$ and, for every $i \in I$, there is some $h_i:M_i \cong N_i$, then $\Pi h_i : \Pi M_i / U \cong \Pi N_i / U$, where $\Pi h_i$ is defined by taking $[i \mapsto f(i)]_U \in \Pi M_i / U$ to $[i \mapsto h_i(f(i))]_U \in \Pi N_i / U$;
	
	\item if $\seq{M_i \in K : i \in I}, \seq{N_i \in K : i \in I}$ and, for every $i \in I$, there is some $h_i : M_i \to N_i$, then $\Pi h_i : \Pi M_i / U \to \Pi N_i / U$, where $\Pi h_i$ is defined by taking $[i \mapsto f(i)]_U \in \Pi M_i / U$ to $[i \mapsto h_i(f(i))]_U \in \Pi N_i / U$; 
	
	\item if $I = \kappa$ and $\seq{M_i \in K : i < \kappa}$ is an increasing sequence, then the ultrapower embedding $h: \bigcup_{i < \kappa} M_i \to \Pi M_i /U$ defined as $h(m) = [f_m]_U$, where
	$$
	f_m(i) = 
	\begin{cases}
	m & \text{if $m \in |M_i|$}\\
	\te{arbitrary} & \text{otherwise}
	\end{cases}
	$$
	is a $K$-embedding; and 
	
	\item if $\seq{M_i \in K : i \in I}$ is a directed set, so in particular $M := \bigcup_{i \in I} M_i \in K$ and, for all $m \in |M|$, we have $[m] = \{ i \in I : m \in M_i\} \in U$, then the ultrapower embedding $h: M \to \Pi M_i / U$ is a $K$-embedding, where $h(m) = [f_m]_U$ and 
	$$
	f_m(i) = 
	\begin{cases}
	m & \text{if $m \in |M_i|$}\\
	\te{arbitrary} & \text{otherwise}
	\end{cases}
	$$

\end{enumerate}
\end{theorem}

{\bf Proof:} If $K$ is an AEC essentially below $\kappa$, then either it is a model of an $L_{\kappa, \omega}$ theory or $LS(K) < \kappa$.  In the first case, this follows from \L o\'{s}' Theorem for $L_{\kappa, \omega}$.

If $LS(K) < \kappa$, then Shelah's Presentation Theorem above says that $K = PC(T_1, \Gamma, L(K))$ for $|T_1| = LS(K) < \kappa$.  During the following proofs, we use the fact observed at \cite{shelahfobook}.VI.0.2 that an ultraproduct of reducts is the reduct of the ultraproducts.
\begin{enumerate}

	\item Each $M_i \in K = PC(T_1, \Gamma, L(K))$, there is some $L(T_1)$ structure $M_i^* \in EC(T_1, \Gamma)$ such that $M_i = M_i^* \rest L(K)$.  Then $\Pi M_i / U = \Pi (M_i^* \rest L(K)) / U = \Pi M_i^* / U \rest L(K)$, such that $\Pi M_i / U$ is the restriction to $L(K)$ of a $L(T_1)$ structure.  Furthermore, there is an $L_{\kappa, \omega}$ sentence $\psi$ st, for any $L(T_1)$ structure $M$, $M \models \psi$ iff $M \in EC(T_1, \Gamma)$.  Thus, for all $i \in I$, $M_i^* \models \psi$.  So by \L o\'{s}' Theorem for $L_{\kappa, \omega}$, $\Pi M^*_i / U \models \psi$.  Thus, $\Pi M^*_i / U \in EC(T_1, \Gamma)$ and $\Pi M_i / U = \Pi M^*_i / U \rest L(K) \in PC(T_1, \Gamma, L(K)) = K$.
	
	\item From  Shelah's Presentation Theorem, for each $N_i$, there are $, M_i^*, N_i^* \in EC(T_1, \Gamma)$ such that $M_i^* \rest L(K) = M_i$, $N_i^* \rest L(K) = N_i$, and $M_i^* \subseteq N_i^*$.  By the above part, $\Pi M_i^* / U, \Pi N_i^* / U \in EC(T_1, \Gamma)$ and, by the definition of an ultraproduct, $\Pi M_i^* / U \subseteq \Pi N_i^* / U$.  Again applying Shelah's Presentation Theorem, we get that
	$$(\Pi M_i^* / U) \rest L(K) = \Pi M_i / U \prec_K \Pi N_i / U = (\Pi N_i^* / U) \rest L(K)$$
	
	\item First we note that $\Pi h_i$ is a bijection.  If $[f]_U \in \Pi N_i /U$, then $[i \mapsto h_i^{-1}(f(i))]_U \in \Pi M_i / U$ and $\Pi h_i([i \to h_i^{-1}(f(i))]_U) = [f]_U$.  If $[f]_U \neq [g]_U \in \Pi M_i / U$, then $\{ i \in I : f(i) = g(i) \} \notin U$.  But this left hand side is $\{ i \in I : h_i(f(i)) = h_i(g(i)) \}$, so $\Pi h_i([f]_U) \neq \Pi h_i ([g]_U)$.\\
	Now we must show that it respects $L(K)$.  Suppose that $R \in L(K)$ is an $n$-ary relation.  Then
	\begin{eqnarray*}
	\Pi M_i / U \models R([f_1]_U, \dots, [f_n]_U) \\
	\{ i \in I : M_i \models R(f_1(i), \dots, f_n(i)) \} \in U \\
	\{ i \in I : N_i \models R(h_i(f_1(i)), \dots, h_i(f_n(i))) \} \in U\\
	\Pi N_i / U \models R(\Pi h_i/U([f_1]_U), \dots, \Pi h_i/U([f_n]_U)) 
	\end{eqnarray*}
	as desired.  The same proof works for functions, or assume $L(K)$ is relational by replacing functions with their graph.
	
	\item For each $i \in I$, we have a $h_i: M_i \cong h_i(M_i)$ with $h_i(M_i) \prec_K N_i$; see the definition of a $K$-embedding.  From above, we know that $\Pi h_i(M_i) / U \prec \Pi N_i / U$ and $\Pi h_i : \Pi M_i / U \cong \Pi h_i (M_i) / U$.  So by the definition of a $K$ embedding, we have our conclusion.
	
	\item This follows from the next one.  Note that, by $\kappa$ completeness, $[m] = \{ \alpha < \kappa : \alpha \geq \beta \} \in U$, where $\beta = \min \{ \gamma < \kappa : m \in |M_\gamma| \}$.
	
	\item Since we modulus by $U$, the definition of $f_m$ only matters on a measure one set, namely $[m]$.  We proceed as in (1) and (2).  We can extend each $M_i$ to $M_i^* \in EC(T_1, \Gamma)$.  Then $\Pi M_i / U = \Pi (M_i^* \rest L) / U = (\Pi M_i^* / U) \rest L$.  Then this induces an $L(T_1)$ expansion of $h(M)$ called $h(M)^* \subseteq \Pi M_i^* \rest L / U$.  So $h: M \to \Pi M_i / U$. \hfill \dag\\

\end{enumerate}

Note, in particular, that in (3) and $(4)$, we have defined the `ultraproduct' of a series of embeddings.  We will generally refer to this as the average of those embeddings and will later use this fact in particular when $N \in K$ and we have many $f_i \in Aut N$; then we know that $\Pi f_i \in Aut \Pi N / U$.

In our definition of essentially below, we hoped to capture all AECs that are closed under complete enough ultraproducts in the sense above.  However, this is not the case: we could take an AEC $K$ which is essentially below $\kappa$ and form the AEC $K' = (K_{\kappa^+})^{up}$ (recall \cite{shelahaecbook}.II.1.23) by taking out all models of size $\kappa$ or smaller.  Then $K'$ is not essentially below $\kappa$, but is still closed under $\kappa$-complete ultraproducts.

However, the hypothesis of an AEC which is essentially below $\kappa$ is natural and somewhat tight, in the sense that there are simple examples of AECs that just fail to be essentially below $\kappa$ and are not closed under $\kappa$-complete ultraproducts.  The following example mirrors the construction of nonstandard models of PA.

\begin{example}
Let $L = \{ <, c_\alpha \}_{\alpha < \kappa}$ and set $\psi \in L_{\kappa^+, \omega}$ to be the sentence
$$\te{``$<$ is a linear order''} \wedge \forall x ( \vee_{\alpha<\kappa} x = c_\alpha) \wedge (\bigwedge_{\alpha < \beta < \kappa} c_\alpha < c_\beta)$$
Let $\eff$ be a $\kappa$-sized fragment containing $\psi$.  Then $K = (\te{Mod } \psi, \prec_\eff)$ is an AEC with $LS(K) = \kappa$, so it `just fails' to be essentially below $\kappa$.  Any $M \in K$ is isomorphic to $(\kappa, \in, \alpha)_{\alpha<\kappa}$.  Thus, $K$ is not closed under $\kappa$-complete ultraproducts.
\end{example}

The results of this paper that have a hypothesis of ``essentially below'' will all continue to hold in any AEC that is closed under sufficiently closed ultraproducts.

Now we are ready to establish the main theorem of this section, that AECs that are essentially below a strongly compact cardinal are tame and type short.  This allows us to connect our large cardinal assumptions to known model theoretic properties.  Afterwards, we will continue our investigation of ultraproducts of AECs; these results will make more sense in light of the fact that types are determined by their $< \kappa$ restrictions.

In this theorem, we assume that $K$ has a monster model.  However, this is not necessary and we do not even need to assume amalgamation for the conclusion.  We include the stronger assumptions to simplify the proof, but provide Theorem \ref{noap} as a ``proof of concept'' that this assumption can be removed.

\begin{theorem}\label{strongcompactnesstheoremtheorem}\label{strongcompactnesstheorem}
Suppose $K$ is essentially below $\kappa$, $\kappa$ is strongly compact, and $K$ has amalgamation, joint embedding, and no maximal models.  Then types are determined by the restrictions of their domain to $< (\kappa + LS(K)^+)$-sized models and their length to $< \kappa$-sized sets; that is, given $M \in K$ and $p, q \in S^I(M)$,
\begin{center}
if $p^{I_0} \rest M_0 = q^{I_0} \rest M_0$ for all $I_0 \in P_\kappa I$ and $M_0 \in P_{\kappa + LS(K)^+}^* M$, then $p = q$.
\end{center}
\end{theorem}

The inclusion of `$+LS(K)^+$' is needed for the case that $K$ is the class of models of some theory in a fragment $\eff$ of $L_{\kappa, \omega}$ with $LS(K) = |\eff| \geq \kappa$; in this case, it would be impossible for $K$ to be $<\kappa$ tame because there would be no models of size $< \kappa$.

First, we prove a technical lemma.  In our proof of Theorem \ref{strongcompactnesstheorem}, there is a place where we will want to take an ultraproduct of our monster model.  However, this would run counter to our intuition of the monster model containing all models since the monster model cannot contain its own ultraproduct.  To avoid this, we introduce a smaller model that functions as the monster model exactly as we need, but without any blanket assumptions of containing all models or being model homogeneous.   We call such a model a \emph{local monster model}.

\begin{lemma}[Local Monster Model] \label{LMM} \label{lmm}
Suppose we have some collection $\{M_i \in K_{\leq \mu} : i < \mu\}$ and $\{f_i \in Aut \sea : i < \mu \}$ such that each $M_i \prec \sea$.  Then there is some $\enn \in K_\mu$ such that for each $i < \mu$ we have $M_i \prec \enn$ and $f_i \rest \enn \in Aut \enn$.
\end{lemma}

\prf Let $N_0 \prec \sea$ of size $\mu$ such that $\bigcup_{i < \mu} |M_i| \subset |N_0|$.  Then each $M_i \prec N_0$.  For $n < \omega$, if we have $N_n$, set $N_{n + 1} \prec \sea$ to be of size $\mu$ such that it contains $|N_n| \cup \bigcup_{i < \mu} (f_i[N_n] \cup f_i^{-1}[N_n])$.  Then set $\enn = \cup_{n<\omega} N_n$. \hfill \dag \\

{\bf Proof of Theorem \ref{strongcompactnesstheoremtheorem}:}  Let $p, q \in S^I(M)$ be as above.  Find $X = \seq{x_i : i \in I} \models p$ and $Y = \seq{y_i : i \in I} \models q$.  Then, by Lemma \ref{LMM}, we find a local monster model $\enn$ such that, for all $(I_0, M_0) \in P_\kappa I \times P_{\kappa+LS(K)^+}^* M$, there is some $f_{(I_0, M_0)} \in Aut_{M_0} \enn$ such that $f_{(I_0, M_0)} (x_i) = y_i$ for all $i \in I_0$.  Next, by the final equivalent definitions of strongly compact cardinals from Definition \ref{strongcompdef}, we find a fine, $\kappa$ complete ultrafilter $U$ on $P_\kappa I \times P_{\kappa+LS(K)^+}^* M$; that is, one such that $[(i, m)] = \{ (I_0, M_0) \in P_\kappa I \times P_{\kappa+LS(K)^+}^* M : i \in I_0, m \in M_0\} \in U$ for all $i \in I$ and $m \in M$.\\
Then, by Theorem \ref{losaec}.6, $\Pi \enn / U \in K$ and our average of these automorphisms $f \in Aut \Pi \enn / U$.  Recall that $f$ takes $[(I_0, M_0) \to g(I_0, M_0)]_U$ to $[(I_0, M_0) \to f_{(I_0, M_0)}(g(I_0, M_0))]_U$.  Now we must prove two claims
\begin{itemize}

	\item {\bf $f$ fixes $h(M)$}\\
	Let $m \in |M|$.  Given any  $i \in I$, $[(i, m)] \in U$ and, if $(I_0, M_0) \in [(i, m)]$, then $m \in M_0$, such $f_{(I_0, M_0)}(m) = m$.  Thus, 
	$$[(i, m)] \subset \{ (I_0, M_0) \in P_\kappa I \times P_{\kappa+LS(K)^+}^* M : f_{(I_0, M_0)}(m) = m \} \in U$$
	and $f \circ h(m) = [(I_0, M_0) \to f_{(I_0, M_0)}(m)]_U = [(I_0, M_0) \to m]_U = h(m)$.
	
	\item {\bf $f(h(x_i)) = h(y_i)$ for every $i \in I$}\\
	Let $i \in I$.  Given any $m \in |M|$, $[(i, m)] \in U$ and, if $(I_0, M_0) \in [(i, m)]$, then $i \in I_0$, so $f_{(I_0, M_0)}(x_i) = y_i$.  Thus, 
	$$[(i, m)] \subset \{ (I_0, M_0) \in P_\kappa I \times P_{\kappa+LS(K)^+}^* M : f_{(I_0, M_0)}(x_i) = y_i \} \in U$$
	and $f \circ h(x_i) = [(I_0, M_0) \to f_{(I_0, M_0)}(x_i)]_U = [(I_0, M_0) \to y_i]_U = h(y_i)$.

\end{itemize}
Now we have the following commutative diagram
\[
 \xymatrix{ \enn  \ar[r]^{f\circ h}  & \Pi \enn / U\\
M \ar[u] \ar[r] & \enn \ar[u]_h
 }
\]
with $f \circ h(x_i) = h(y_i)$ for all $i \in I$.  Thus, $p = q$. \hfill \dag\\

The above theorem can be interpreted as saying that if we have two different types, then they are different on a ``formula,'' if we take formula to mean a type of $< \kappa$ length over a domain of size $< \kappa$.  With this definition of formula, we can replace a large type by the set consisting of all of its small restrictions and type equality will be preserved.  In the rest of this section, we will see that, since $\kappa$ is strongly compact, this notion of formulas as small types will be fruitful.  We now return to the development of our ultraproducts with a version of \L o\'{s}' Theorem.

Also, we strengthen our hypothesis to $LS(K) < \kappa$ instead of just $K$ essentially below $\kappa$.  This is because \cite{makkaishelah}.2.10 has shown that, with a monster model, Galois types in models of a $L_{\kappa, \omega}$ theory correspond to consistent sets of formulas from a fragment of $L_{\kappa, \kappa}$, so the following results are already known.

Note that the following theorem only requires a measurable cardinal.

\begin{theorem}[\L o\'{s}' Theorem for AECs, part 2]
Suppose that $K$ is an AEC with amalgamation, joint embedding, and no maximal models and $\kappa$ is a measurable cardinal such that $\kappa < LS(K)$.  Let $N^- \prec N \in K$ and $p \in S(N^-)$ with $\|N^-\| + \ell(p) < \kappa$ and $U$ be a $\kappa$-complete ultrafilter on $I$.  Then $[g]_U \in \Pi N/U$ realizes $h(p)$ iff $\{ i \in I :g(i) \vDash p\} \in U$, where $h: N \to \Pi N/U$ is the canonical ultrapower embedding.
\end{theorem}

{\bf Proof:} \underline{$\Leftarrow$}  Suppose that $[g]_U \in \Pi N/U$ with $X := \{ i \in I : g(i) \vDash p \} \in U$. Let $a \vDash p$.  By Lemma \ref{lmm}, there is a local monster model $\enn$ such that, for each $i \in X$, there is $f_i \in Aut_{N^-} \enn$ such that $f_i(g(i)) = a$.

Define $f^+ :\Pi \enn/U \to \Pi \enn/U$ to be the average of these maps.  That is, $f^+$ of $[i \mapsto k(i)]_U$ is $[i \mapsto f_i(k(i))]_U$; although these $f_i$'s don't exist everywhere, they exist on a $U$-large set and this is enough.  Then, by \L o\'{s}' Theorem for AECs, $f^+ \in Aut_{h(N^-)} \Pi \enn/U$.  Note that $h(N^-)$ is $h'' N^-$ and not $\Pi N^- / U$.  Also, $f_i$ sends $g(i)$ to $a$ on a large set, so $f^+([g]_U) = h(a)$.  So $[g]_U$ realizes $tp(h(a)/h(N^-)) = h(p)$ as desired.

\underline{$\Rightarrow$}  Let $[g]_U \in \Pi N / U$ realize $h(p) \in S(h(N^-))$.  For each $q \in S(N^-)$, set $X_q = \{ i \in I : g(i) \textrm{ realizes } q \}$.  Since different $q$'s are mutually exclusive, these are all disjoint and they partition $I$.  We easily have $|S(N^-)| \leq 2^{\|N^-\|} < \kappa$.  Since $U$ is $\kappa$-complete, this means that, for some $q_0 \in S(N^-)$, $X_{q_0} \in U$.  By the previous direction, that means that $[g]_U$ realizes $h(q_0)$.  But by assumption, the type of $[g]_U$ over $h(N^-)$ is $h(p)$, so $p = q_0$.  Thus $X_p = \{ i \in I : g(i) \textrm{ realizes } p \} \in U$, as desired. \hfill \dag\\

Now that we have \L o\'{s}' Theorem, we prove a companion result to Theorem \ref{strongcompactnesstheorem}.  This motivated our conception of types as sets of smaller types or ``formulas.''  Here we show that, as with the first-order case, any consistent set of formulas can be completed to a type, even when the set is incomplete.

We introduce some notation to make this as general as possible.  Even with a strongly compact cardinal, a key difference between small types and formulas is that there is no negation of a type: given a type $p$, there is no type $q$ such that all elements realize either $p$ or $q$ and not both.  To compensate for this, we want to allow specification of both types to be realized and types to be avoided.  In the following, $X$ represents the types to be realized and $\neg X$ represents the types to be avoided.

\begin{defin}
Fix $M \in K$ and $I$ a linear order.  Let $X \subset \{ p \in S^{I_0}(M^-) : I_0 \in P_\kappa I\te{ and }M^- \in P^*_\kappa M\}$ and $\neg X \subset \{ q \in S^{I_0}(M^-) : I_0 \in P_\kappa I\te{ and }M^- \in P^*_\kappa M\}$.
\begin{itemize}
	\item We say that $\ba = \seq{a_i : i \in I}$ realizes $(X, \neg X)$, written $\ba \vDash (X, \neg X)$ iff, for every $p \in X$, $\seq{a_i : i \in \ell(p)} \vDash p$ and, for every $q \in \neg X$, $\seq{a_i : i \in \ell(q)} \not\vDash q$.  We say that $(X, \neg X)$ is consistent iff it has a realization.
	\item We say that $(X, \neg X)$ is $< \kappa$ consistent iff $(X_0, \neg X_0)$ is consistent for every $X_0 \in P_\kappa X$ and $\neg X_0 \in P_\kappa \neg X$.
\end{itemize}
\end{defin}

\begin{theorem}
Suppose $K$ is an AEC with amalgamation, joint emebedding, and no maximal models and $\kappa$ is strongly compact such that $LS(K) < \kappa$.  Let $M\in K$ and $I$ be a linear order.  Given $X \subset \{ p \in S^{I_0}(M^-) : I_0 \in P_\kappa I\te{ and }M^- \in P^*_\kappa M\}$ and $\neg X \subset \{ q \in S^{I_0}(M^-) : I_0 \in P_\kappa I\te{ and }M^- \in P^*_\kappa M\}$, $(X, \neg X)$ is consistent iff it is $< \kappa$ consistent.
\end{theorem}

{\bf Proof:} One direction is obvious, so suppose that $(X, \neg X)$ is $<\kappa$ consistent.  For every $N \in P^*_\kappa M$, let $X_N = \{ p \in X : \dom p \prec N \}$ and $\neg X_N = \{ q \in \neg X : \dom q \prec N \}$.  Then, by assumption $(X_N, \neg X_N)$ is consistent, so there is $\ba_N = \seq{a^N_i : i \in I}$ that realizes $(X_N, \neg X_N)$; if $X_N = \neg X_N = \emptyset$, then pick $\ba_N$ arbitrarily.  Let $M^+ \succ M$ contain all $\ba_N$ and let $U$ be a $\kappa$ complete, fine ultrafilter on $P^*_\kappa M$.  Recall that $h: M^+ \to \Pi M^+ / U$ is the canonical embedding.

For each $i \in I$, set $a_i := [N \mapsto a_i^N]_U$ for $i \in I$ and set $\ba := \seq{a_i : i \in I}$.  We claim that $\ba \vDash h(X, \neg X) = (\{h(p) : p \in X\}, \{h(q) : q \in \neg X\})$.  Suppose $p \in X$ and $M^- \prec M$ and $I_0 \subset I$ such that $p \in S^{I_0}(M^-)$.  Then
$$[M^-] = \{ N \in P^*_\kappa M : M^- \prec M \} \subset \{ N \in P^*_\kappa M : \seq{a_i^N : i \in I_0} \vDash p \}$$
by construction.  Since this first set is in $U$ by fineness, $\seq{a_i : i \in I_0} \vDash h(p)$ by the previous theorem.  Now suppose $q \in \neg X$ and $M^- \prec M$ and $I_0 \subset I$ such that $q \in S^{I_0}(M^-)$.  For contradiction, suppose that $\seq{a_i : i \in I_0} \vDash h(q)$.  Then, by the previous theorem,  $\{ N \in P^*_\kappa M : \seq{a^N_i : i \in I_0} \vDash q\} \in U$.  Then let $N' \in \{ N \in P^*_\kappa M : \seq{a^N_i : i \in I_0} \vDash q\}  \cap [M^-]$; this intersection is nonempty because it is in $U$.  Then $\seq{a^{N'}_i : i \in I_0}$ both realizes and does not realize $q$, a contradiction.  Thus, $\seq{a_i : i \in I_0}$ does not realize $h(q)$ and we have shown $\ba \vDash h(X, \neg X)$, as desired.

Let $h^+$ be an $L(K)$ isomorphism that extends $h$ and has image $\Pi M^+ / U$.  Then $(h^+)^{-1}(\ba)$ witnesses the consistency of $(X, \neg X)$. \hfill \dag\\

We can use a similar argument to transfer saturation from $M$ to $\Pi M / U$.

\begin{theorem}
Suppose $K$ is an AEC with amalgamation, joint embedding, and no maximal models and $\kappa$ is strongly compact such that $LS(K) < \kappa$.  For all $M \in K$ and linear order $I$, there is some $\kappa$ complete $U$ such that, for any $p \in S^I(M)$ that has all $< \kappa$ restrictions realized in $M$, $\Pi M / U \models h(p)$.
\end{theorem}

{\bf Proof:}  Let $U$ be a $\kappa$ complete, fine ultrafilter on $P_\kappa I \times P^*_\kappa M$ and, for each small approximation $p^{I_0} \rest M^-$, pick some $\ba^{(I_0, M^-)}:=\seq{a^{(I_0, M^-)}_i \in M : i \in I_0} \vDash p^{I_0} \rest M^-$.  Now consider the sequence $\seq{[(I_0, M_0) \mapsto a^{(I_0, M^-)}_i]_U : i \in I}$.  This sequence is in $\Pi M / U$ since each $a^{(I_0, M_0)}_i \in M$.  By the same argument as the previous theorem, this sequence realizes $h(p)$. \hfill \dag\\

\begin{cor}
If $M \in K$ is $< \kappa$ saturated, then there is some $\kappa$-complete $U$ such $\Pi M /U$ realizes all types over $M$.
\end{cor}

\section{Measurable} \label{measurablesection}

We now turn our attention to what happens if our large cardinal is only measurable.

\begin{defin}[\cite{jech}.17]
An uncountable cardinal $\kappa$ is \emph{measurable} iff there is a normal, $\kappa$-complete ultrafilter on $\kappa$.  Equivalently, there is some elementary embedding $j: V \to \emm$ with critical point $\kappa$ such that ${}^\kappa \emm \subset \emm$.
\end{defin}

Unsurprisingly, we don't get as strong results here.  Instead, we just get results of $(<\lambda, \lambda)$-tameness and  type shortness whenever $\cf \lambda = \kappa$.  Reexamining the above proof, an argument readily presents itself by using the $\kappa$-complete ultrafilter on $\kappa$ and redoing the above arguments.  Instead of repeating the above proof, we prove this theorem in two different ways: once with a monster model and using the ultrapower definition, and the second time using ultrafilters but no assumption of amalgamation at all.  We do these proofs in order to showcase different large cardinal techniques on AECs.  The use of the elementary embedding is of particular interest, because this is the formulation of large cardinals most studied by modern set theorists and will hopefully shed light on future work in this direction, while the proof without amalgamation shows the we get the results from just large cardinals and do not need additional, structural assumptions on $K$, like amalgamation.\\

\begin{theorem}
Suppose $K$ is an AEC essentially below $\kappa$ measurable with amalgamation, joint emebedding, and no maximal models.  Let $M = \bigcup_{\alpha < \kappa} M_\alpha$ and $I = \bigcup_{\alpha < \kappa} I_\alpha$ and $p \neq q \in S^I(M)$.  Then, there is some $\alpha_0 < \kappa$ such that $p^{I_{\alpha_0}} \rest M_{\alpha_0} \neq q^{I_{\alpha_0}} \rest M_{\alpha_0}$.
\end{theorem}

{\bf Proof:} Let $M = \bigcup_{\alpha < \kappa} M_i$ and $p \neq q \in S(M)$, as above.  Let $X = \seq{x_i : i \in I}$ and $Y = \seq{y_i : i \in I}$ realize $p$ and $q$ respectively.  Since $\kappa$ is measurable, there is some normal, $\kappa$-complete ultrafilter $U$ on $\kappa$ such that we get the following commuting and elementary diagram
\[
\xymatrix{ V \ar[rr]^j \ar[dr]_i & & \emm \\
& \Pi V/ U \ar[ur]_\pi &
}
\]
where $i$ is the ultrapower embedding, $\pi$ is the Mostowski collapse, $crit j = \kappa$, and ${}^\kappa \emm \subset \emm$.  Since $V$ is the set-theoretic universe, we also have $\emm, \Pi V / U \subset V$.  Similarly, $j(I) = \bigcup_{\alpha < j(\kappa)} I_\alpha'$ and $I_\kappa' = \bigcup_{\alpha < \kappa} j(I_\alpha)$.  Set $X_\alpha' = \seq{x_i : i \in I_\alpha'}$ and $Y'_\alpha = \seq{y_i : i \in I_\alpha'}$.\\
By elementarity, we have that $j(p) \neq j(q) \in (S^{j(I)}(j(M)))^\emm$ and $j(M) = \bigcup_{\alpha < j(\kappa)} M_\alpha'$, where $\seq{M_\alpha' : \alpha < j(\kappa)} = j(\seq{M_\alpha : \alpha < \kappa})$ and, for $\alpha < \kappa$, $M_\alpha' = j(M_\alpha)$.  $\seq{M_\alpha' : \alpha < j(\kappa)}$ is continuous, so $M_\kappa' = \bigcup_{\alpha < \kappa} M_\alpha' = \bigcup_{\alpha < \kappa} j(M_\alpha)$.\\
A priori, all that is known is that $j(M) \in j(K)$, which $\emm$ thinks is an AEC.  In fact, $j(K) = K^\emm$; this follows from Theorem \ref{losaec} since, for any $N \in K$, $i(N) = \Pi N / U \in K$.  Since $\pi$ is an isomorphism, $\pi \circ i(N) = j(N) \in K$.\\
For any $N \in K$, we note that $i``N \in K$ is isomorphic to $N$ and has universe $\{ [\alpha \to n]_U : n \in |N|\}$, so $i``N \prec i(N) = \Pi N / U$ by the above.  So $i``M_\alpha \prec i(M_\alpha)$ for every $\alpha < \kappa$.  Thus $i``M_\alpha \prec i(M_\beta)$ for every $\beta \geq \alpha$ and, taking a union over the $\beta < \kappa$, $i``M_\alpha \prec \bigcup_{\beta < \kappa} i(M_\beta)$.  Now taking a union over $\alpha < \kappa$, $\bigcup_{\alpha < \kappa} i``M_\alpha = i``M \prec \bigcup_{\beta < \kappa} i(M_\beta)$.  Applying $\pi$ to both sides yields 
$$j``M = (\pi \circ i)``M \prec \pi( \bigcup_{\beta < \kappa} i(M_\beta)) \prec \bigcup_{\beta<\kappa} \pi \circ i(M_\beta) = M_\kappa ' $$
Since $M_\kappa' \prec j(M)$, which is the domain for $j(p)$ and $j(q)$, we have that $j(p) \rest j''M$ and $j(q) \rest j''M$ are defined.  Similarly, $j``I \subset I_\kappa'$ so $j``X \subset X_\kappa'$ and $j``Y \subset Y_\kappa'$.\\
We wish to show $j(p)^{I_\kappa'} \rest M_\kappa'$ and $j(q)^{I_\kappa'} \rest M_\kappa'$ are different.  So we compute
\begin{eqnarray*}
i(p) &=& tp(i(X)/i(M)\\
&=& tp((\Pi X /U) / (\Pi M / U))
\end{eqnarray*}
and $i(q) = tp((\Pi Y/U) / (\Pi M / U))$.  However, we know that $tp(X/M) \neq tp(Y/M)$.  So $tp(i``X/i``M) \neq tp(i``Y/i``M)$.  Since $i``M \prec \Pi M / U$, $i``I \subset \Pi I/U$ and non-equality of types goes up, we have that $i(p)^{i``I} \rest i``M \neq i(q)^{i``I} \rest i``M$.  Applying our isomorphism $\pi$, we get
\begin{eqnarray*}
\pi(i(p)^{i``I} \rest i``M) &\neq& \pi(i(q)^{i``I} \rest i``M)\\
\pi \circ i (p)^{(\pi \circ i)``I} \rest (\pi \circ i)``M &\neq& \pi \circ i (q)^{(\pi \circ i)``I} \rest (\pi \circ i)``M \\
j(p)^{j``I} \rest j``M &\neq& j(q)^{j``I} \rest j``M
\end{eqnarray*}
as desired.  Since $j``M \prec M_\kappa'$ and $j``I \subset I_\kappa'$, we have that $j(p)^{I_\kappa'} \rest M_\kappa' \neq j(q)^{I_\kappa'} \rest M_\kappa'$.\\
So far, we have argued completely in $V$.  However, since equality of types is existentially witnessed and a witness in $\emm$ would also be a witness in $V$, this holds true in $\emm$ as well.  So, we get the following
\begin{eqnarray*}
\emm &\models& j(p)^{I_\kappa'} \rest M_\kappa' \neq j(q)^{I_\kappa'} \rest M_\kappa'\\
\emm &\models& \exists \alpha < j(\kappa) \te{ st for $N =$ the $\alpha$th member of } j(\seq{M_\beta : \beta < \kappa}) \te{ and $J =$ the $\alpha$th member of } \\
& & j(\seq{I_\beta : \beta < \kappa}), j(p)^J \rest N \neq j(q)^J \rest N\\
V &\models& \exists \alpha < \kappa \te{ st for $N =$ the $\alpha$th member of } \seq{M_\beta : \beta < \kappa}  \te{ and $J =$ the $\alpha$th member of } \\
& & \seq{I_\beta : \beta < \kappa}, p^J \rest N \neq q^J \rest N\\
V &\models& \exists \alpha < \kappa \te{ st} p^{I_\alpha} \rest M_\alpha \neq q^{I_\alpha} \rest M_\alpha
\end{eqnarray*}
Since $V$ is the universe, there is some $\alpha_0 < \kappa$ such that $p^{I_{\alpha_0}} \rest M_{\alpha_0} \neq q^{I_{\alpha_0}} \rest M_{\alpha_0}$. \hfill \dag\\

\begin{cor}
$K$ is fully $(< \lambda, \lambda)$-tame and fully $(<\lambda, \lambda)$-type short whenever $cf \lambda = \kappa$ and $\lambda > LS(K)$.
\end{cor}

Finally, we wish to weaken the assumptions on the theorems above to remove the use of the monster model.  Note that, because, in these contexts, we can always take an ultrapower and $M \precneqq \Pi M / U$ for any $M \in K$ at least the size of the completeness of the ultrafilter, we already have no maximal models.  So, in particular, we remove the assumptions of amalgamation and joint embedding.  The loss of amalgamation is particularly worrisome because it is used to prove that $\sim_{AT}$ is an equivalence relation and we only have that $\sim$ is a \emph{non-trivial} transitive closure of $\sim_{AT}$.  Also, we now use the complete strength of the closure theorem for ultraproducts of AECs.

\begin{theorem} \label{noap}
Suppose $K$ is an AEC essentially below $\kappa$ measurable.  $K$ is fully $(< \lambda, \lambda)$ tame and fully $(<\lambda, \lambda)$ type short for $\lambda > LS(K)$ with $\cf \lambda = \kappa$.
\end{theorem}

{\bf Proof:}  For ease, we only show tameness.  Type shortness follows similarly, but would add extra notation to an already notation heavy proof.
Let $M \in K_\lambda$ and let $p, q \in S(M)$ such that $p \rest N = q \rest N$ for all $N \in P^*_\lambda M$.  Find $a, b, N^0, N^1$ such that $p = tp(a/M, N^0)$ and $q = tp(b/M, N^1)$.  Then we can find resolutions $\seq{M_i, N_i^0, N_i^1 \in K_{< \lambda} : i < \kappa}$ of $M$, $N^0$, and $N^1$ respectively such that $a \in |N_0^0|$ and $b \in |N_0^1|$.  Then we know that, for all $i < \kappa$, $(a, M_i, N_i^0) \sim (b, M_i, N_i^1)$.  Since $\sim$ is the transitive closure of $\sim_{AT}$, for every $i < \kappa$, there is some $n_i < \omega$ such that, for all $\ell \leq n_i$, we have $N_i^{\frac{\ell}{n_i}}$ and $a_\ell^i$ such that $a_0^i = a$, $a_{n_i}^i = b$, every $a_\ell^i \in  |N_i^{\frac{\ell}{n_i}}|$, and 

$$(a, M_i, N_i^0) \sim_{AT} (a^i_1, M_i, N_i^\frac{1}{n_i}) \sim_{AT} \cdots \sim_{AT} (a^i_{n_i-1}, M_i, N_i^\frac{n_i - 1}{n_i}) \sim_{AT} (b, M_i, N_i^1)$$

Since there are only countably many choices for $n_i$ and $\cf \kappa > \omega$, there is some $n$ that occurs cofinally often; WLOG, we may thin our sequence and assume $n_i = n$ for all $i < \kappa$.  In particular, note that Theorem \ref{losaec}.5 does not require continuity.
Now, by the definition of $\sim_{AT}$, for all $i < \kappa$ and $\ell < n$, there is some $N_{i, l}^* \succ N_i^\frac{\ell + 1}{n}$ and $f_{i, \ell}: N_i^\frac{\ell}{n} \to N^*_{i, \ell}$ such that

\[
 \xymatrix{\ar @{} [dr] N_i^\frac{\ell}{n}  \ar[r]^{f_{i, \ell}}  & N^*_{i, \ell}\\
M_i \ar[u] \ar[r] & N_i^\frac{\ell+1}{n} \ar[u]
 }
\]

commutes and $f_{i, \ell}(a^i_\ell) = a^i_{\ell + 1}$.  Looking across all $\ell < n$, we get the following commuting diagram

\[
\xymatrix{
 & N^*_{i, 1} & \cdots & N_{i, n - 2}^* & \\
 N^*_{i, 0} & N_i^\frac{1}{n} \ar[l] \ar[u]_{f_{i, 1}} & \cdots & N_i^\frac{n-1}{n} \ar[u] \ar[r]^{f_{i, n -1}} & N^*_{i, n - 1} \\
 & N_i^0 \ar[ul]_{f_{i, 0}} & M_i \ar[r] \ar[ur] \ar[ul] \ar[l] & N_i^1 \ar[ur] & 
}
\]
so $f_{i, \ell}(a^i_\ell) = a^i_{\ell+1}$ for all $i < \kappa$ and $\ell < n$.

Let $U$ be some $\kappa$-complete ultrafilter on $\kappa$.  Now we take the ultraproduct of the above diagrams.  Recall that we use $\Pi f_i$ to denote the average of maps $f_i$; see the discussion after the proof of Theorem \ref{losaec}.
\[
\xymatrix{
 & \Pi N^*_{i, 1} / U & \cdots & \Pi N_{i, n - 2}^* / U & \\
 \Pi N^*_{i, 0} / U & \Pi N_i^\frac{1}{n} / U \ar[l] \ar[u]_{\Pi f_{i, 1} } & \cdots & \Pi N_i^\frac{n-1}{n} / U \ar[u] \ar[r]^{\Pi f_{i, n -1}} & \Pi N^*_{i, n - 1} / U\\
 & \Pi N_i^0 / U \ar[ul]^{\Pi f_{i, 0} } & \Pi M_i / U \ar[r] \ar[ur] \ar[ul] \ar[l] & \Pi N_i^1 / U \ar[ur] & 
}
\]
Also by our hypotheses, if we take the function $h: M \to \Pi M_i / U$ given by $h(m) = [ i \to m]_U$, then this is a $K$ embedding.  Note that, although the function $i \to m$ is not well-defined for all $i$, by $U$'s $\kappa$ completeness, it is defined on a measure one set, so the $h$ is still well-defined.  We can similarly define $h_0: N_0 \to \Pi N_i^0 / U$ and $h_1: N_1 \to \Pi N_i^1 / U$.  Note that, for all $m \in M$, $h(m) = h_0(m) = h_1(m)$.  These allow us to construct the following commutative diagram
\[
\xymatrix{
 & \Pi N^*_{i, 1} / U & \cdots & \Pi N_{i, n- 2}^* / U & \\
 \Pi N^*_{i, 0} / U & \Pi N_i^\frac{1}{n} / U \ar[l] \ar[u]_{\Pi f_{i, 1} } & \cdots & \Pi N_i^\frac{n-1}{n} / U \ar[u] \ar[r]^{\Pi f_{i, n -1}} & \Pi N^*_{i, n - 1} / U\\
  \Pi N_i^0 / U \ar[u]^{\Pi f_{i, 0} } & & \Pi M_i / U \ar[ul] \ar[ur]  & & \Pi N_i^1 / U \ar[u]  \\
 & N^0 \ar[ul]^{h_0} & & N^1 \ar[ur]_{h_1} & \\
 & & M \ar[ul] \ar[ur] \ar[uu]_{h} & &
} 
\]
This is essentially the diagram that we want, but we have to do some renaming to get it into the desired form.  For each $1 \leq \ell < n$, set $\ba_\ell = [ i \mapsto a_\ell^i]_U$ and $\bar N^*_\ell \prec \Pi N_i^{\frac{\ell}{n}} / U$ of size $\lambda + \ell(a)$ containing $\ba_\ell$ and $h(M)$.  Then find some $L(K)$ isomorphism $f_\ell$ that contains $h$ with range $\bar N^*_\ell$ and $N^*_\ell$ such that $f_\ell : N^*_\ell \cong \bar N^*_\ell$.  Set $a_\ell = f_\ell^{-1}(\ba_\ell) \in |N_\ell^*|$.  This gives us the diagram

\[
\xymatrix{
 & \Pi N^*_{i, 1} / U & \cdots & \Pi N_{i, n - 2}^*/U & \\
\Pi N^*_{i, 0} / U & N_1^* \ar[l]_{f_1} \ar[u]_{\Pi f_{i, 1} \circ f_1} & \cdots & N^*_{n-1} \ar[u]_{f_{n-1}} \ar[r]^{\Pi f_{i, n -1} \circ f_{n-1}} & N^*_{i, n - 1} \\
 & N^0 \ar[ul]_{\Pi f_{i, 0} \circ h_0} & M \ar[r] \ar[ur] \ar[ul] \ar[l] & N^1 \ar[ur]_{h_1} & 
}
\]

This witnesses that

$$(a, M, N^0) \sim_{AT} (a_1, M, N_1^*) \sim_{AT} \cdots \sim_{AT} (a_{n-1}, M, N_{n - 1}^*) \sim_{AT} (b, M, N^1)$$

So $(a, M, N^0) \sim (b, M, N^1)$ and $p = q$.\hfill \dag\\

\section{Weakly Compact}

In this section, we establish a number of downward reflection principles using indescribable cardinals.  Tameness follows because it is a downward reflection of type inequality, but these principles apply to many other AEC properties as well.

\begin{defin}[Indescribable Cardinals, \cite{kanamori}.1.6]
\begin{enumerate}

	\item[]
	
	\item For $m, n < \omega$, a cardinal $\kappa$ is \emph{$\Pi^m_n$-indescribable} iff for any $R \subset V_\kappa$ and $\Pi_n^m$-statement $\phi$ in the language of $\{ \in, R \}$, if $\seq{V_\kappa, \in, R} \vDash \phi$, then there is $\alpha < \kappa$ such that $\seq{V_\alpha, \in, R \cap V_\alpha} \vDash \phi$
	
	\item $\kappa$ is \emph{totally indescribable} iff $\kappa$ is $\Pi^m_n$-indescribable for all $n, m < \omega$.

\end{enumerate}	
\end{defin}

Although the indescribability definition is stated in terms of a single $R \subset V_\kappa$, a simple coding argument shows that it is equivalent to allow finitely many $R_0, \dots, R_n \subset V_\kappa$ in the expanded language.

\begin{remark}[\cite{kanamori}]
An uncountable cardinal $\kappa$ is weakly compact iff $\kappa$ is $\Pi^1_1$ indescribable.  Another definition is that any $\kappa$ sized set of sentences from $L_{\kappa, \kappa}$ is consistent iff all of its $< \kappa$ sized subsets are.  For context, if $\kappa$ is measurable then it is $\Pi^2_1$ indescribable and, moreover, for any normal ultrafilter $U$ on $\kappa$, $\{\alpha< \kappa : \alpha \te{ is totally indescribeable} \} \in U$.
\end{remark}

In the following lemma, we are going to code models of an AEC $K$ with $LS(K) < \kappa$ as a subset of $V_\kappa$.  In order to do this, we use the fact that there are two definable functions $g$ and $h$ so
\begin{itemize}

	\item $g: \kappa \to P_\omega \kappa$ is a bijection such that, for all $\mu < \kappa$, we have that $g \rest \mu:\mu \to P_\omega \mu$ is a bijection; and
	
	\item $h: \kappa \times {}^{LS(K)} 2 \to V_\kappa$ is an injection such that, for all $2^{LS(K)} < \mu < \kappa$, we have that $h \rest \mu : \mu \times {}^{LS(K)} 2 \to V_\mu$ is an injection.

\end{itemize}

\begin{lemma}[Coding Lemma]
Suppose $K$ is an AEC such that $LS(K) < \kappa$.  There is $C_K \subset V_\kappa$ and $\Pi^0_n$ formulas $\phi(x), \psi(x, y), \sigma(x, y), \tau(x, y, z), \tau^+(x, y, z) \in L(\{ \in, C_K\})$ such that for $\alpha \leq \kappa$ and $X, Y, f \subset V_\alpha$ and $a \in V_\alpha$, we have
\begin{itemize}
	
	\item $\seq{V_\alpha, \in, C_K \cap V_\alpha} \vDash \phi(X) \iff$\\
$C_K$ decodes from $X$ an $L(K)$ structure $M_X$ and $M_X \in K_{|\alpha|}$
	
	\item $\seq{V_\alpha, \in, C_K \cap V_\alpha} \vDash \psi(X, Y) \iff$\\
$C_K$ decodes from $X$ and $Y$ $L(K)$ structures $M_X$ and $M_Y$, $M_X, M_Y \in K_{|\alpha|}$, and $M_X \prec_K M_Y$

	\item $\seq{V_\alpha, \in, C_K \cap V_\alpha} \vDash \sigma(X, a) \iff$\\
$C_K$ decodes from $X$ $L(K)$ structures $M_X \in K_{|\alpha|}$ and $a \in M_X$

	\item $\seq{V_\alpha, \in, C_K \cap V_\alpha} \vDash \tau(X, Y, f) \iff$\\
$C_K$ decodes from $X$ and $Y$ $L(K)$ structures $M_X$ and $M_Y$ such that $M_X, M_Y \in K_{|\alpha|}$, and $f: M_X \to M_Y$

	\item $\seq{V_\alpha, \in, C_K \cap V_\alpha} \vDash \tau^+(X, Y, f) \iff$\\
$C_K$ decodes from $X$ and $Y$ $L(K)$ structures $M_X$ and $M_Y$, $M_X, M_Y \in K_{|\alpha|}$, and $f: M_X \cong M_Y$

\end{itemize}
\end{lemma}

{\bf Proof:} In the statement, we make reference to a ``decoded structure,'' which we will explain.  By Shelah's Presentation Theorem, we know $K = PC(T_1, \Gamma, L(K))$.  Additionally, we can code $\prec_K$ as an AEC with L\"{o}wenheim-Skolem number $LS(K)$: set $K_{\prec} = \{ (M, |M_0|) : M_0 \prec_K M \}$ and $\prec_{K_\prec} = \{ ((M, |M_0|), (N, |N_0|) \in K_\prec \times K_\prec :M \prec_K N$ and $ M_0 \prec_K N_0 \}$.  Then we have that $K_\prec = PC(T_2, \Gamma', L(K)')$.  WLOG, we can assume that these objects are in $V_{(2^{LS(K)})^{+\omega}}$ and $L(T_1) = \seq{R_i : i < |L(T_1)|}$ and $L(T_2) = \seq{S_i : i < |L(T_2) }$  are relational.  Set $C_K = ((2^{LS(K)})^{+\omega}, L(K), \Gamma, T_1, \Gamma', T_2)$.  Define a $\Pi^0_n$ formula $\phi$ such that $\seq{V_\kappa, \in, C_K \cap V_\alpha} \models \phi(X)$ asserts all of the following
\begin{enumerate}

	\item $C_K$ is an ordered sextuple whose first element is an ordinal; this guarantees that the $V_\alpha$ that models it is above $2^{LS(K)}$ and, thus, can see the other elements.
	
	\item $X$ is in the range of $h$ and $(h^{-1})''X$ is of the form $\{ (i, f_i) : i < \alpha\}$.  Set $C_i = \{ j \in \alpha : f_j(i) = 1\}$.
	
	\item $g''C_0$ should be a set of singletons; denote $\bigcup g''C_0$ by $|M_X|$.
	
	\item $g''C_i$ should be a set of tuples whose length match the arity of $R_i$; denote this set $R^{M^+_X}_i$.
	
	\item $M^+_X = \seq{|M_X|, R^{M^+_X}_i}_{i < LS(K)}$ models $T_1$ and omits each $p \in \Gamma$.
	
	\item Finally, $M_X$ is the model $M_X^+ \rest L(K)$.

\end{enumerate}

Thus, $\seq{V_\kappa, \in, C_K \cap V_\alpha} \models \phi(X)$ iff $M_X \in K$ by Shelah's Presentation Theorem.

For $\psi(X, Y)$, we do a similar decoding process with $T_2$ and $\Gamma'$.

For $\sigma(X, a)$, we need to say that $a$ is in the image of our decoding of $C_0$, which requires a quantifier over an element of $X$.

For $\tau^+(X, Y, f)$, we use $\phi$ to determine that $X$ and $Y$ are codes for elements of our $PC$ class and then say that $f$ is an isomorphism, which again just quantifies over elements of our models and $L$, all of which we have given.

For $\tau(X, Y, f)$, we have a definable way to talk about the image of $X$ under $f$ and combine $\psi$ and $\tau^+$ to say that $f$ is an isomorphism between $X$ and its image and that $X$'s image is a $\prec_K$ submodel of $Y$.\hfill \dag\\

Now we are ready to begin proving theorems from this coding.

\begin{theorem}[Tameness Down for $\Pi^1_1$]
Suppose $K$ is an AEC such that $LS(K) < \kappa$ with $\kappa$-AP and $\kappa$ being $\Pi^1_1$-indescribable.  Then $K$ is $(< \kappa, \kappa)$-tame for $< \kappa$-types.
\end{theorem}

\prf Let $C_K$ be as in the Coding Lemma.  Let $M \in K_\kappa$ and $p \neq q \in S(M)$.  Then we have $p = tp(a/M, N_1)$ and $q = tp(b/M, N_2)$ for $M \prec N_1, N_2 \in K_\kappa$ and $a \in N_1$ and $b \in N_2$.  WLOG, $|N_1| \cup |N_2| \subset V_\kappa$.  Then
\begin{eqnarray*}
\models & & M, N_1, N_2 \in K_\kappa\\
& \wedge& M \prec N_1, N_2 \\
&\wedge &a \in N_1, b \in N_2 \wedge \forall N^* \in K_\kappa, \forall f_i:N_1 \to N^* 
 ( \forall m \in M (f_1(m) = f_2(m)) \to f_1(a) \neq f_2(b) )
\end{eqnarray*}

Let $X, Y_1, Y_2 \subset V_\kappa$ code $M, N_1, N_2$, respectively, according to $C_K$.  Then we rewrite the above as

\begin{eqnarray*}
 \seq{V_\kappa , \in, C_K, X, Y_1, Y_2, \{a\}, \{b\}} &\models& \phi(X) \wedge \phi(Y_1) \wedge \phi(Y_2) \wedge \psi(X, Y_1) \wedge \psi(X, Y_2) \wedge \\
 &\wedge& \sigma(X_1, a) \wedge \sigma(X_2, b) \\
 &\wedge& \forall Y^*, f_1, f_2 \subset V_\kappa [(\phi(Y^*) \wedge \tau(Y_1, Y^*, f_1) \wedge \tau(Y_2, Y^*, f_2) \wedge \\
 & &\wedge [ \forall x \in V_\kappa \sigma(X, x) \to (f_1(x) = f_2(x))] ) \to (f_1(a) \neq f_2(b))]
\end{eqnarray*}

Since everything is first-order except for the single universal quantifier over subsets of $V_\kappa$, this is a $\Pi^1_1$ statement.  So it reflects down to some $\alpha < \kappa$.  Since for this to happen, $\{a\} \cap V_\alpha$ and $\{b\} \cap V_\alpha$ must be nonempty, we must have $a, b < \alpha$.  Let $X' = X\cap V_\alpha$, $Y_1' = Y_1 \cap V_\alpha$, and $Y_2' = Y_2 \cap V_\alpha$.  Then we have that $tp(a/M_{X'}, N_1') = p \rest M_{X'}$ and $tp(b/M_{X'}, N_2') = q \rest M_{X'}$.\\
{\bf Claim:} $p \rest M_{X'} \neq q \rest M_{X'}$\\
If not, then there is some $N^* \in K_{|\alpha|}$ and $f_i:N_i' \to N^*$ that witnesses this with $f_1(m) = f_2(m)$ for all $m \in M$ and $f_1(a) = f_2(b)$.  However, WLOG, $|N^*| \subset \alpha$, so we can code $N^*$ as $Y^* \subset V^\alpha$ according to $C_K$.  Then $f_1, f_2 \subset V_\alpha$ and $Y^*, f_1, f_2$ serve as a counterexample for our downward reflection.\\
Thus, we have our $M_{X'} \in K_{< \kappa}$ such that $p$ and $q$ differ on their restriction to $M_{X'}$. \hfill \dag\\

Above, we assumed amalgamation to simplify the exposition.  However, we could drop this assumption without difficulty by adding a (first-order) quantifier to see how many steps it might take to show $p$ and $q$ are equal.\\

A similar argument gives us a result for type shortness.

\begin{theorem}[Tameness Down for $\Pi^1_1$]
Suppose $K$ is an AEC such that $LS(K) < \kappa$ with $\kappa$-AP and $\kappa$ being $\Pi^1_1$-indescribable.  Then $K$ is $(< \kappa, \kappa)$-type short over $< \kappa$-sized models.
\end{theorem}

This method is not just useful for tameness and type shortness.  It can be used to reflect many AEC properties down.  Only the amount of indescribibility required changes from property to property.  For instance,

\begin{theorem}[Unbounded Categoricity Down for $\Pi^1_2$]
Suppose $K$ is an AEC such that $LS(K) < \kappa$ with $\kappa$ being $\Pi^1_2$-indescribable.  Then for every $\lambda < \kappa$, there is some $\lambda < \mu < \kappa$ such that $K$ is $\mu$-categorical.
\end{theorem}

\prf Let $\lambda < \kappa$. Code $K$ by $C_K$.  We want to find $\lambda < \mu < \kappa$ such that that $K$ is $\mu$-categorical.  Since $K$ is $\kappa$ categorical,
\begin{eqnarray*}
 &\models& \forall M, N \in K_\kappa, \exists f:M \cong N\\
\seq{V_\kappa, \in, C_K, \{\lambda^+\}} &\models& \forall X, Y \subset V_\kappa \exists f \subset V_\kappa [\phi(X) \wedge \phi(Y) \to \tau^+(X, Y, f)] \wedge \exists x (x \in \{\lambda^+\})
\end{eqnarray*}
Then this reflects down to some $\alpha < \kappa$.  Since $V_\alpha \cap \{ \lambda^+ \}$ is not empty, we get that $\alpha > \lambda^+$, so $|\alpha| > \lambda$.  Set $\mu = |\alpha|$ and let $M, N \in K_\mu$.  WLOG, $|M|, |N| \subset \alpha$, so we can code these by $X$ and $Y$, respectively.  Then
$$\seq{V_\alpha, \in, C_K \cap V_\alpha} \models \phi(X) \wedge \phi(Y)$$
Since our statement of categoricity reflects down to $\alpha$, there is some $f \in V_\alpha$ such that $f:M \cong N$.  So $K$ is $\mu$ categorical.\hfill \dag\\

\begin{remark}
Recalling what has been said about work on Shelah's Categoricity Conjecture, one may initially hope that this downward reflection might be massaged to make the downward reflection hold at a successor cardinal.  However, this is unlikely, since successor (and singular limit) cardinals are necessarily first-order describable, so all we could guarantee of $\mu$ is that it is strongly inaccessible.
\end{remark}

We have many other theorems of this type:

\begin{theorem}[Unbounded (Disjoint) Amalgamation Down for $\Pi^1_2$]
Suppose $K$ is an AEC such that $LS(K) < \kappa$ with $\kappa$ (disjoint) amalgamation and $\kappa$ being $\Pi^1_2$-indescribable.Then, for every $\lambda < \kappa$, there is some $\lambda < \mu < \kappa$ such that $K$ has the $\mu$ (disjoint) amalgamation.
\end{theorem}

\begin{theorem}[Unbounded Uniqueness of Limit Models Down for $\Pi^2_1$] \label{limitmodels}
Suppose $K$ is an AEC such that $LS(K) < \kappa$ with $\kappa$ being $\Pi^2_1$-indescribable. If $K_\kappa$ has a unique limit model, then, for every $\lambda < \kappa$, there is some $\lambda < \mu < \kappa$ such that $K_\mu$ has a unique limit model.
\end{theorem}

The general heuristic for determining how much indescribability is required to transfer a property of an AEC down is to look at the quantifiers needed to state this property and translate quantifiers over elements to $\Pi^0$ quantifiers; over models or embeddings to $\Pi^1$ quantifiers; and over sequences of models or embeddings to $\Pi^2$ quantifiers.  Following this, sequences of sequences of models would require $\Pi^3$ quantifiers, but there seem to be no useful AEC properties requiring a quantifier of this sort.

\section{Conclusion} \label{conclusion}

cy of Shelah's Eventual Categoricity Conjecture for Successors by combining our results with those of \cite{tamenessthree} and \cite{sh394}.  After doing so, we apply our results to other results in the literature.

Before we can apply the results of of \cite{tamenessthree} and \cite{sh394}, we must show that categoricity implies their hypotheses of no maximal models, joint embedding, and amalgamation.  If $K$ is the class of models of some $L_{\kappa, \omega}$ sentence, then this is done in \cite{makkaishelah}.\S1.  We generalize these arguments to an AEC $K$ with $LS(K) < \kappa$ by introducing the notion of universal closure as a generalization of existential closure.

\begin{defin}
$M \in K$ is called universally closed iff given any $N \prec M$ and $N' \succ N$, both of size less than $\kappa$, if there is $M^+ \succ M$ and $g: N' \to_N M^+$, then there is $f: N' \to_N M$.
\end{defin}

We omit the parameter $\kappa$ from the name because it will always be fixed and clear from context.  Note that if there is an $M^+$ witnessing that $M$ is not universally closed, then there is one of size $\|M\|$.

Recall that $M$ is an amalgamation base when all $M_1$ and $M_2$ extending $M$ can be amalgamated over $M$.

\begin{lemma}\label{uc=ab}
Suppose $K$ is an AEC and $\kappa$ is strongly compact such that $LS(K) < \kappa$.  Then any universally closed $M \in K_{\geq \kappa}$ is an amalgamation base.
\end{lemma}

{\bf Proof:}  Let $M$ be universally closed and $M \prec M_1, M_2$.  First, we show we can amalgamate every small approximation of this system.  Let $N \prec M$ and $N_\ell \prec M_\ell$ such that $N \prec N_\ell$ for $\ell = 1, 2$ with $N, N_1, N_2 \in K_{< \kappa}$.  Then $M_\ell$ is an extension of $M$ such that $N_\ell$ can be embedded into it over $N$.  Since $M$ is universally closed, there is $f_\ell:N_\ell \to_N M$.  Find $N_* \prec M$ of size $< \kappa$ such that $f_1(N_1), f_2(N_2) \prec N_*$.  Then this is an amalgamation of $N_1$ and $N_2$ over $N$.

Now we will use our strongly compact cardinal.  Set
$$X = \{ \bN = (N^\bN, N_1^\bN, N_2^\bN) \in (K_{< \kappa})^3 : N^\bN \prec M, N_\ell^\bN \prec M_\ell, N^\bN \prec N_\ell^\bN \te{ for } \ell = 1, 2\}$$
For each $\bN \in X$, the above paragraphs shows that there is an amalgam of this triple.  Fix $f^\bN_\ell: N^\bN_\ell \to N_*^\bN$ to witness this fact.  For each $(A, B, C) \in [M]^{<\kappa} \times [M_1]^{< \kappa} \times [M_2]^{<\kappa}$, define
$$[(A, B, C)]:= \{ \bN \in X:A \subset N^\bN, B \subset N_1^\bN, C \subset N^\bN_2\}$$
These sets generate a $\kappa$-complete filter on $X$, so it can be extended to a $\kappa$-complete ultrafilter $U$.  By \L o\'{s}' Theorem for AECs, since this ultrafilter is fine, we know that the ultrapower map $h$ is a $K$-embedding, so 
\begin{eqnarray*}
h:M \to \Pi N^\bN/U & &\\
h_\ell:M_\ell \to \Pi N_\ell^\bN/U & & \te{ for }\ell = 1, 2
\end{eqnarray*}
Since these maps have a uniform definition, they agree on their common domain $M$.  Furthermore, we can average the $f^\bN_\ell$ maps to get
$$\Pi f^\bN_\ell: \Pi N_\ell^N / U \to \Pi N_*^\bN / U$$
and the maps agree on $\Pi N^\bN/U$ since each of the individual functions do.  Then we can put these maps together to get the following commutative diagram that witnesses the amalgamation of $M_1$ and $M_2$ over $M$.

\[
 \xymatrix{
 & \Pi N^\bN_2/U \ar[rr]_{\Pi f_2^\bN} & & \Pi N_*^\bN/U \\
M_2 \ar[ur]_{h_2} & & &  \\
 & \Pi N^\bN/U \ar[rr] \ar[uu] & & \Pi N^\bN_1/U \ar[uu]_{\Pi f_1^\bN}\\
M \ar[uu] \ar[ur]_h \ar[rr] & & M_1\ar[ur]_{h_1} &  
  }
\]
\hfill \dag\\

Now we use this result to derive the needed properties from categoricity.  We focus on the case where $K$ is categorical in $\lambda$ of cofinality at least $\kappa$ because it is simpler and suffices for our application.  However, the methods of \cite{makkaishelah} can extend these results to categoricity in other cardinals.  We use here the result of Solovay \cite{solovay} that $\cf \mu \geq \kappa$ implies $\mu^{< \kappa} = \mu$ when $\kappa$ is strongly compact.

\begin{prop} \label{getthatmonster}
Suppose $K$ is an AEC such that $LS(K) < \kappa$ strongly compact.  If $K$ is categorical in $\lambda$ such that $\cf \lambda \geq \kappa$, then $K_{\geq \kappa}$ has amalgamation, joint embedding, and no maximal models.
\end{prop}

{\bf Proof:}  $K_{\geq \kappa}$ has no maximal models by \L o\'{s}' Theorem for AECs, since a model can be strictly embedded into its ultraproduct.  This doesn't use categoricity and only needs $\kappa$ to be measurable.

For joint mapping, we can use categoricity and no maximal models to get joint mapping below and at the categoricity cardinal.  Above the categoricity cardinal, we use amalgamation and categoriciy.  This relies only on the other properties and not directly on any large cardinals.

For amalgamation, we use the above result that universally closed models are amalgamation bases.

First, we show that a universally closed model exists in any cardinal $\mu$ of cofinality at least $\kappa$, which includes the categoricity cardinal.  Let $M \in K_\mu$ and consider all possible isomorphism types of $N \prec N'$ from $K_{< \kappa}$ with $N \prec M$.  There are at most $\mu^{< \kappa} \cdot 2^{< \kappa} = \mu$ many such types.  We enumerate them $(N_\alpha, N'_\alpha)$ for $\alpha < \mu$.  Set $M = M_0$.  Then for each $\alpha < \mu$, if there is some $M_\alpha^+ \succ M_\alpha$ of size $\mu$ such that there is $g: N_\alpha' \to_{N_\alpha} M_\alpha^+$ but no $f: N_\alpha' \to_{N_\alpha} M_\alpha$, then set $M_{\alpha+1} = M_\alpha^+$.  Otherwise, $M_{\alpha+1} = M_\alpha$.  At limit $\alpha$, we take limits of the increasing chain.  Set $M^* = \cup_{\alpha < \mu} M_\alpha \in K_\mu$.

Now we iterate this process $\kappa$ many times: set $M^0 = M$, $M^{\alpha + 1} = (M^\alpha)^*$, and $M^\alpha = \cup_{i < \beta} M^i$ for limit $\alpha \leq \kappa$.  Then, $M^\kappa$ is universally closed.  By $\lambda$ categoricity, this means that every model in $K_\lambda$ is universally closed.

Second, we show that every model in $K_{> \lambda}$ is a universally closed.  Let $M \in K_{> \lambda}$.  Suppose that there are $N\prec  N' \in K_{< \kappa}$ and $M^+ \succ M \succ N$ and $g: N' \to_N M^+$.  Let $M' \prec M$ be of size $\lambda$ and contain $N$.  Then, by the above, $M'$ is universally closed with $M^+ \succ M'$, so there is $f: N' \to_N M'$.  Then $f: N' \to_N M$.  Since $N$ and $N'$ were arbitrary, $M$ is universally closed.

Third, we show that all models in $K_{\geq \kappa}$ are amalgamation bases and, thus, $K_{\geq \kappa}$ has the amalgamation property.  Let $M \prec M_1, M_2$.  If $M \in K_{\geq \lambda}$, then $M$ is universally closed and, thus, an amalgamation base by Lemma \ref{uc=ab}.  If not, then we can find some $\kappa$ complete ultrafilter $U$ and take an ultraproduct to get a proper extension
\[
\xymatrix{
& \Pi M_2 / U & & \\
M_2 \ar[ur]_h & & &\\
 & \Pi M / U \ar[rr] \ar[uu]  & & \Pi M_1 / U \\
M \ar[rr] \ar[uu] \ar[ur]_h& & M_1 \ar[ur]_h &
}
\]

This is a larger triple of models that, if we could amalgamate it, would give us an amalgamation of $M_1, M_2$ over $M$.  Then, we can continue to take ultrapowers of this triple, taking direct limits at unions, until the base model has size at least $\lambda$.  Then, by the above, it must be an amalgamation base, so we can amalgamate $M_1$ and $M_2$ over $M$.

Thus, all models in $K_{\geq \kappa}$ are amalgamation bases, so $K_{\geq \kappa}$ has the amalgamation property. \hfill \dag\\

Now that we have amalgamation, joint embedding, and no maximal models, we can generalize the result of \cite{makkaishelah} to all AECs essentially below a strongly compact.

\begin{theorem} \label{bigresult}
Suppose $\kappa$ is a strongly compact cardinal and $K$ is an AEC essentially below $\kappa$.  If $K$ is categorical in some successor $\lambda^+$ greater than $\kappa^+ + LS(K)^+$, then it is categorical in all $\mu \geq \min \{ \lambda^+, \beth_{(2^{Hanf(LS(K))})^+} \}$.
\end{theorem}

{\bf Proof:} By Theorem \ref{strongcompactnesstheorem}, $K$ is $< (\kappa + LS(K)^+)$ tame, so it is $\kappa+LS(K)^+$ tame.  Then, $K_{\geq \kappa}$ is an AEC with $LS(K_{\geq \kappa}) = \kappa$ that is $\kappa$-tame.  Additionally, by Proposition \ref{getthatmonster}, $K$ has amalgamation, joint embedding, and no maximal models.  Thus, by \cite{tamenessthree}.5.2, we know that $K$ is categorical for every $\mu \geq \lambda^+$.  Then $K$ is definitely categorical in a successor above $\beth_{(2^{Hanf(LS(K))})^+}$.  So, by \cite{sh394}.9.5, it is categorical everywhere down to $\beth_{(2^{Hanf(LS(K))})^+}$. \hfill \dag\\

Note that the downward categoricity transfer result from \cite{sh394} does not use any tameness assumption.  This result shows that given an AEC wtih amalgamation that is categorical in a successor cardinal $\lambda$ above $\beth_{(2^{Hanf(LS(K))})^+}$, this AEC is also categorical in all cardinals in the interval $[\beth_{(2^{Hanf(LS(K))})^+}, \lambda]$.

Now we show that Shelah's Eventual Categoricity Conjecture for Successors follows from large cardinal assumptions:

\begin{theorem} \label{shelahcatconj}
If there are proper class many strongly compact cardinals, then Shelah's Eventual Categoricity Conjecture for Successors holds.
\end{theorem}

{\bf Proof:} Let $\lambda$ be a cardinal and pick $\mu_\lambda = \min \{ \mu^+ : \mu \geq \lambda$ and $\mu$ is strongly compact $\}$  Note that $\beth_{(2^{Hanf(\lambda)})^+} < \mu_\lambda$.  If $K$ is categorical in some successor $\mu$ above $\mu_\lambda$, then Theorem \ref{bigresult} implies that $K$ is categorical everywhere above $\mu_\lambda$. \hfill \dag\\

While the hypothesis of this theorem seems very strong, we do note that  \cite{jech}.20.22 and .24 show that the consistency of it follows from the existence of an extendible cardinal $\lambda$; in fact, $V_\lambda$ is a model of the hypothesis.

Beyond the categoricity result, \cite{makkaishelah} introduces a very well behaved independence relation similar to the first-order notion of coheir.  This nonforking notion is generalized to AECs in Boney and Grossberg \cite{shorttamedep} and its uniqueness is established in Boney, Grossberg, Kolesnikov, and Vasey \cite{bgkv}.  Of particular note is that no large cardinal hypothesis is need, only the conclusions of Theorem \ref{strongcompactnesstheorem} for a specific AEC.

Of particular interest in the proof of \ref{shelahcatconj} is that we get, from the hypothesis of a proper class of strongly compact cardinals, the conclusion that \emph{every} AEC with arbitrarily large models is eventually tame.  Examining the ZFC counterexamples of \cite{hash323} \cite{untame}, the proven failure of tameness occurs at some small level bounded by $\aleph_\omega$.  However, these classes have arbitrarily large models, so our results can apply.  In particular, if there is a strongly compact cardinal, these classes exhibit the strange behavior of being $(\aleph_0, \aleph_k)$-tame very low, failing to be $(\aleph_k, \aleph_{k+1})$-tame, and then becoming $<\kappa$ tame at the strongly compact cardinal. 

Turning to measurable cardinals, \cite{kosh} derive amalgamation from categoricity and \cite{measure2}  proves a downward categoricity transfer in $L_{\kappa, \omega}$.  However, the papers do not use the specifics of $L_{\kappa, \omega}$ beyond that it is closed under $\kappa$ complete ultralimits, see \cite{kosh}.1.7.1.  The methods of Theorem \ref{losaec} can be used to show closure under these ultralimits as well.  Thus, we can extend their work to get the following results:

\begin{theorem}
Suppose $K$ is an AEC such that $LS(K) < \kappa$ measurable.  If $K$ is categorical in some $\lambda \geq \kappa$, then
\begin{enumerate}

	\item $K_{[LS(K) + \kappa, \lambda)} = \{ M \in K : LS(K) + \kappa \leq \|M\| < \lambda \}$ has the amalgamation property; and
	
	\item if $\lambda$ is also a successor above $\beth_{(2^{LS(K)})^+}$, then $K$ is categorical in all $\mu$ with $\beth_{(2^{LS(K)})^+} \leq \mu \leq \lambda$.

\end{enumerate}
\end{theorem}

Beyond ultralimits, stronger large cardinals have more complicated constructions that witness their existence, such as extenders for strong cardinals \cite{jech}.20.28.  Again, arguments similar to Theorem \ref{losaec} will show closure under these constructions as well for AECs essentially below them.

In Theorem \ref{limitmodels}, we mention limit models.  While not discussed more in this paper, these are well-studied objects and the uniqueness of limit models seems to be an important dividing line for AECs; see \cite{gvv}, \cite{vandierennomax} \cite{nomaxerrata}, or \cite{shvi635} for more information.

\section{Further work} \label{furtherwork}
\label{futureworksection}

As always, new answers lead to new questions.

In this paper, we have shown that the following statements follow from different large cardinals:

\begin{enumerate}
 \item[$(*)_\kappa^-$] Every AEC $K$ with $LS(K) < \kappa$ is $(<\kappa, \kappa)$-tame.
 \item[$(*)_\kappa$] Every AEC $K$ with $LS(K) < \kappa$ is $<\kappa$-tame.
 \item[$(*)$] Every AEC $K$ with arbitrarily large models is $<\lambda$-tame in some $\lambda > LS(K)$.
\end{enumerate}

We proved the same results for type shortness, but we focus this discussion on tameness because more is known.

A natural investigation is into these properties on their own.  Can they hold at small cardinal?  If so, do they have large cardinal strength?

A basic first result is that none of these properties can hold at $\aleph_k$ for $k < \omega$.  This follows from the Hart-Shelah examples \cite{hash323} \cite{untame}.

A second result is that $(*)_\kappa^-$ for $\kappa$ regular and not weakly compact implies $V \neq L$.  To see this, first recall that Baldwin and Shelah \cite{nonlocality} construct an AEC that is not $(< \kappa, \kappa)$ tame from an almost free, non-free, non-Whitehead group of size $\kappa$.  In $L$, such a group is known to exists at precisely the non-weakly compact, regular cardinals; see Eklof and Mekler's book \cite{ekloffmekler}.  Combining these two facts, we have our proof.  The construction in Eklof and Mekler has two main steps:
\begin{itemize}
	\item non-reflecting stationary sets are used to construct almost free, non-free groups of every cardinality; and
	\item weak diamond on every stationary set is used to inductively show that all Whitehead groups are free.
\end{itemize}
The non-reflecting stationary sets suggest a natural tension with the compactness of the cardinals used in this paper.  However, also being non-Whitehead seems to be a crucial part of the Baldwin and Shelah construction.  It is not currently known if non-tameness follows just from almost free, non-freeness nor what the precise conditions for the existence of an almost free, non-free, non Whitehead group are.

A potential first step in achieving $(*)_\kappa^-$ and the other properties at small cardinals is the work of Magidor and Shelah in \cite{almostfree}.  Starting from $\omega$ many supercompact cardinals, they construct
\begin{enumerate}
 \item a model where every $\aleph_{\omega^2+1}$-free group is $\aleph_{\omega^2+2}$-free; and
 \item a model where every $\kappa$-free group is free for $\kappa = \min \{ \lambda \in \te{CARD} : \lambda = \aleph_\lambda\}$.
\end{enumerate}
In the first model, there is no known candidate for a counterexample to $(*)_{\aleph_{\omega^2+2}}^-$ and, in the second, there is no known candidate for a counterexample to $(*)_{\aleph_\kappa}$.  Further investigation will be needed to determine if these properties hold or if there are more non-tame AECs.

\end{document}